\pdfoutput=1
\documentclass{amsart}
\usepackage{xcolor}
\usepackage{hyperref}
\hypersetup{
	colorlinks=true,
	citecolor=blue,
	linkcolor=blue,
	filecolor=magenta,      
	urlcolor=cyan,
}
\usepackage{zref-clever}
\zcsetup{cap,noabbrev,nameinlink=false}

\zcsetup{
  countertype={
    thm=theorem,
    cthm=theorem,
    ccor=corollary,
    cor=corollary,
    que=question,
    lem=lemma,
    prop=proposition
  }
}

\AddToHook{env/thm/begin}{\zcsetup{countertype={equation=theorem}}}
\AddToHook{env/thm/end}{\zcsetup{countertype={equation=}}}
\AddToHook{env/prop/begin}{\zcsetup{countertype={equation=proposition}}}
\AddToHook{env/prop/end}{\zcsetup{countertype={equation=}}}
\AddToHook{env/lem/begin}{\zcsetup{countertype={equation=lemma}}}
\AddToHook{env/lem/end}{\zcsetup{countertype={equation=}}}
\AddToHook{env/cor/begin}{\zcsetup{countertype={equation=corollary}}}
\AddToHook{env/cor/end}{\zcsetup{countertype={equation=}}}
\AddToHook{env/dfn/begin}{\zcsetup{countertype={equation=definition}}}
\AddToHook{env/dfn/end}{\zcsetup{countertype={equation=}}}
\AddToHook{env/rem/begin}{\zcsetup{countertype={equation=remark}}}
\AddToHook{env/rem/end}{\zcsetup{countertype={equation=}}}
\AddToHook{env/exa/begin}{\zcsetup{countertype={equation=example}}}
\AddToHook{env/exa/end}{\zcsetup{countertype={equation=}}}
\AddToHook{env/conj/begin}{\zcsetup{countertype={equation=conjecture}}}
\AddToHook{env/conj/end}{\zcsetup{countertype={equation=}}}

\usepackage[all, knot]{xy}\SelectTips{eu}{}

\usepackage{amsmath,amscd}
\usepackage{amssymb}
\usepackage{mathtools}
\usepackage{stmaryrd}
\usepackage{url}
\usepackage{tikz-cd}
\usepackage{enumitem}

\usepackage{setspace}
\usepackage{microtype}

\usepackage{times}

\usepackage{tikz}

\newtheorem{thm}[equation]{Theorem}
\newtheorem*{thm*}{Theorem}
\newtheorem*{cor*}{Corollary}
\newtheorem*{dfn*}{Definition}

\newtheorem{cthm}{Theorem}

\newtheorem{que}[cthm]{Question}

\newtheorem{cor}[equation]{Corollary}
\newtheorem{prop}[equation]{Proposition}
\newtheorem{lem}[equation]{Lemma}

\theoremstyle{definition}

\newtheorem{rem}[equation]{Remark}

\newcommand{\inj}{\hookrightarrow}

\newcommand{\opn}{\operatorname}
\newcommand{\cat}[1]{\operatorname{\mathsf{#1}}}

\newcommand{\mfrak}[1]{\mathfrak{#1}}

\newcommand{\mrm}[1]{\mathrm{#1}}

\renewcommand{\k}{\Bbbk}

\newcommand{\m}{\mfrak{m}}

\newcommand{\p}{\mfrak{p}}

\newcommand{\injdim}{\operatorname{inj\,dim}}
\newcommand{\projdim}{\operatorname{proj\,dim}}
\newcommand{\flatdim}{\operatorname{flat\,dim}}

\def\skewtimes{\ltimes\!}

\newcommand{\restl}{{}_R \opn{Rest}}
\newcommand{\restr}{\opn{Rest}_R}

\newcommand{\restlb}{{}_B \opn{Rest}}
\newcommand{\restrb}{\opn{Rest}_B}

\title{Gorenstein symmetry can fail for noetherian rings}

\author{Liran Shaul}

\address{Department of Mathematics, University of Haifa, Haifa 3498838, Israel}

\email{shaul@math.haifa.ac.il}


\numberwithin{equation}{section}

\begin{document}

\begin{abstract}
We construct a two-sided noetherian ring $A$ which has injective dimension $1$
as a left module over itself and infinite injective dimension as a right module
over itself. This shows that Gorenstein symmetry can fail for two-sided
noetherian rings. We show that such an example exists over any field of
characteristic different from $2$, that it has a dualizing complex, and that its
derived category of right modules does not satisfy injective generation in the
sense of Rickard-Keller. Finally, we show that our example satisfies the
Auslander condition, giving a negative answer to a question of Huang.
\end{abstract}

\maketitle

\setcounter{tocdepth}{1}
\tableofcontents

\setcounter{section}{-1}
\section{Introduction}

In this paper, all rings are associative and unital,
and modules are by default left modules. 
Thus, an $A$-module is a left $A$-module,
and a right $A$-module is an $A^{\opn{op}}$-module.

In commutative algebra contexts, a commutative noetherian ring $A$ of finite Krull dimension is called Gorenstein if it has finite injective dimension over itself. The study of the Gorenstein property encompasses an important subject in modern commutative algebra and algebraic geometry.
An important class of Gorenstein rings in commutative algebra are the regular rings. A commutative noetherian ring of finite Krull dimension is regular if it has finite global dimension.

When passing to noncommutative contexts, a symmetry question arises:

\begin{que}
Let $A$ be a ring, and suppose that $\injdim_A( \phantom{ }_A A) < \infty$.
Does it follow that $\injdim_A (A_A) < \infty$?
\end{que}

Similarly one can ask this question for global dimension.
The global dimension question has a well known classical answer:
there are examples of rings which have finite global dimension on the left and infinite global dimension on the right. See \cite{Jat}. 
However, when restricting attention to the class of \textbf{two-sided noetherian} rings, Auslander showed \cite[Corollary 5]{Aus} that the left and the right global dimension coincide. In particular, if one of them is finite, then so is the other one.

Going back to symmetry of injective dimension, the main focus of this paper, Zaks showed in \cite[Lemma A]{Zaks} that if $A$ is a two-sided noetherian ring, and if both $\injdim_A( \phantom{ }_A A) < \infty$ and 
$\injdim_A (A_A) < \infty$, then
$\injdim_A( \phantom{ }_A A) =  \injdim_A (A_A)$.
In that case, $A$ is called two-sided Gorenstein, 
or two-sided Iwanaga-Gorenstein, following the foundational work \cite{IW1,IW2} of Iwanaga about this class of rings.

In \cite[Section 2]{KKS}, the authors state they know of no example of a two-sided noetherian ring for which this symmetry fails. The purpose of this paper is to construct such an example.
Thus, our main contribution is to show that in the noncommutative noetherian world, left-Gorenstein and two-sided Gorenstein are not the same. Here by left-Gorenstein we mean a two-sided noetherian ring which has finite injective dimension as a left module over itself.

Our main result says:
\begin{cthm}\label{cthm:main}
There exists a two-sided noetherian ring $A$ with the property that
\[
\injdim_A( \phantom{ }_A A) = 1, \quad \injdim_A(A_A) = +\infty.
\]
\end{cthm}

For artin algebras, the question if injective dimension on the left and on the right coincide is called the Gorenstein symmetry conjecture,
see \cite{AusRei,BeKr,BReit}. For artin algebras,
and more generally, artinian rings, the conclusion of \zcref{cthm:main} cannot hold, because by a result of Hoshino \cite{Ho},
artinian rings with injective dimension at most $1$ on the left are two-sided Gorenstein. By \cite[Theorem 15.1]{Lam2}, a two-sided noetherian ring having injective dimension $0$ on one side has it on both sides.
Thus, the value $1$ in the above result is the smallest possible for such a counterexample.

Let us describe the contents of the paper in more detail.
The ring $A$ in \zcref{cthm:main} will be a twisted ring of dual numbers, $A = R[x;\sigma]/(x^2)$ over some ring $R$. 
The ring $R$ will be required to be a hereditary ring,
equipped with some injective ring homomorphism $\sigma:R \to R$ which is not surjective. In \zcref{sec:gen}, we consider such an $R$ and such a $\sigma$, which moreover satisfies that
\begin{equation}\label{eqn:condInt}
1 \in (R:\sigma(R))_r\cdot R,
\end{equation}
where 
\[
(R:\sigma(R))_r = \{x \in R \mid R\cdot x \subseteq \sigma(R)\},
\]
is the right conductor of the inclusion of rings $\sigma(R) \subseteq R$, 
a notion recalled in \zcref{sec:gen}.
We then show in \zcref{thm:main-gen-const} that $R[x;\sigma]/(x^2)$ satisfies the conclusions of \zcref{cthm:main}.
Note that if $R$ were commutative, then \zcref{eqn:condInt} would imply that $\sigma$ is surjective. Thus, the main technical heart of this paper is that one can construct such an $R$ and such a $\sigma:R \to R$ which satisfy \zcref{eqn:condInt}, while $\sigma$ is not surjective.

The ring $R$ is a subring of a polynomial ring $D[t]$ over a certain division ring $D$, which will also be required to satisfy various technical conditions. The construction of $D$ is carried out in \zcref{sec:div}, and the exact conditions we need $D$ to satisfy are spelled in \zcref{thm:divisionRing}.

With $D$ in hand, we move to constructing $R$ and $\sigma$. 
This is done in \zcref{sec:R}, whose main result, \zcref{thm:mainR}
says that over any field $\k$ with $\opn{char}(\k) \ne 2$,
there is a two-sided noetherian hereditary $\k$-algebra $R$,
together with an injective $\k$-algebra homomorphism $\sigma:R \to R$ which is not surjective, which satisfies \zcref{eqn:condInt}.

After that, \zcref{sec:dual} is dedicated to showing that the algebra $A$ has a dualizing complex, in the sense of Yekutieli \cite{Yek}. This is shown in \zcref{thm:dualcomplex}.
Combining all of these, we obtain the following refinement of \zcref{cthm:main}:

\begin{cthm}\label{cthm:main3}
Let $\k$ be any field with $\opn{char}(\k) \ne 2$.
Then there exists a $\k$-algebra $A$ which satisfies the following properties:
\begin{enumerate}
\item $A$ is a two-sided noetherian ring.
\item $A$ is not two-sided artinian.
\item $A$ has a dualizing complex.
\item It holds that:
\[
\injdim_A( \phantom{ }_A A) = 1, \quad \injdim_A(A_A) = +\infty.
\]
\end{enumerate}
\end{cthm}

The construction of the dualizing complex is needed in order for $A$ to satisfy the covariant Grothendieck duality \cite{InKr,Neeflat},
which in turn implies that our recent characterization \cite{ShGor} of Gorenstein symmetry in terms of the behavior of the derived category of $A$ holds. This is discussed in detail in \zcref{sec:RK}.
There, we study questions about injective generation and projective cogeneration in the derived categories associated to $A$.
These notions are recalled in that section, and its main result,
\zcref{thm:failure} states:

\begin{cthm}
Let $\k$ be any field with $\opn{char}(\k) \ne 2$,
and let $A$ be the $\k$-algebra constructed in \zcref{cthm:main3}.
Then the following additional properties hold for $A$:
\begin{enumerate}
\item Injectives generate for $A$ on the left.
\item Injectives do not generate for $A$ on the right.
That is, $\opn{Loc}_{\cat{D}(\opn{Mod}_A)}(\opn{Inj}_A) \ne \cat{D}(\opn{Mod}_A)$.
\item $\opn{Coloc}_{\cat{D}({}_A \opn{Mod})}({}_A \opn{Proj}) \ne \cat{D}({}_A \opn{Mod})$. 
That is, projectives do not cogenerate for $A$ on the left.
\item $\opn{Coloc}_{\cat{D}({}_A \opn{Mod})}({}_A \opn{Flat}) \ne \cat{D}({}_A \opn{Mod})$. 
That is, flats do not cogenerate for $A$ on the left.
\item There exists a bounded above cochain complex $J$ of injective left $A$-modules which is acyclic but not contractible.
\item There exists a bounded below cochain complex $P$ of projective right $A$-modules which is acyclic but not contractible.
\item There exists an injective left $A$-module $I$ with infinite flat dimension.
\item There exists a dualizing complex $\omega'$ over $A$ which has infinite flat dimension as a complex of left $A$-modules.
\end{enumerate}
\end{cthm}

Finally, in \zcref{sec:Auslander}, we recall the Auslander condition. For many classes of rings, it is known that a ring which satisfies the Auslander condition and has finite injective dimension on one side must have finite injective dimension on both sides. These classes include all artin algebras, all artinian rings, and all two-sided noetherian rings which are module-finite over a commutative noetherian center. The main result, \zcref{cor:Auslander}, states:

\begin{cthm}
Let $\k$ be any field with $\opn{char}(\k) \ne 2$,
and let $A$ be the $\k$-algebra constructed in \zcref{cthm:main3}.
Then $A$ satisfies the Auslander condition.    
\end{cthm}

This provides a negative answer to a question of Huang \cite{Huang},
which asked if the Auslander condition implies Gorenstein symmetry for two-sided noetherian rings.

\section{A general construction violating Gorenstein symmetry}\label{sec:gen}

Throughout this section, $R$ is an arbitrary associative unital ring,
and $\sigma:R \to R$ is a ring homomorphism (in particular, $\sigma(1) = 1$).
We will further assume that $\sigma$ is injective but not surjective.
Specific instances of $R$ and $\sigma$ will be constructed in the next section,
but for now we work in this generality.

Following \cite[Section 1.2.3]{NCBook},
we consider the skew polynomial ring $B = R[x;\sigma]$.
Let us recall the construction.
An element of $B$ has a unique expression as a finite sum
\[
\sum_{k=0}^n x^k a_k, \quad a_k \in R.
\]
and the multiplication satisfies $a\cdot x = x\cdot \sigma(a)$ for all $a \in R$.
Let 
\begin{equation}\label{eqn:hero}
A = B/(x^2) = R[x;\sigma]/(x^2).
\end{equation}

The ring $A$ will be the main hero of this section.

We can give an alternative description of the ring $A$ in terms of trivial extensions. Let $M$ be the $R$-$R$-bimodule whose underlying set is equal to $R$, and with $R$-actions given by 
\[
r\circ m := \sigma(r) \cdot m,  \quad r\in R, m \in M
\]
and
\[
m \circ r := m\cdot r \quad m \in M, r \in R
\]
We will sometimes denote $M$ by $\phantom{ }_{\sigma} R$.
We then consider the trivial extension ring $\widetilde{A} = R \skewtimes M$, in the sense of \cite{TExtBook}. 
Recall, by definition that as a set, $\widetilde{A} = R\times R$.
We shall denote the element $(0,1) \in \widetilde{A}$ by $\epsilon$.
Thus, an element of $\widetilde{A}$ has the form $x_1 + \epsilon y_1$, with $x_1,y_1 \in R$,
$M = \epsilon\cdot R$,
and the multiplication on $\widetilde{A}$ is given by
\begin{equation}\label{eqn:multA}
(x_1 + \epsilon y_1) \cdot (x_2 + \epsilon y_2) = x_1\cdot x_2 + \epsilon(\sigma(x_1)\cdot y_2+ y_1\cdot x_2).
\end{equation}
There is a ring isomorphism $A \to \widetilde{A}$,
given by $x\mapsto \epsilon$ and restricting to the identity on $R$. We will use it to identify these two constructions,
and, by abuse of notation, will denote both of them by $A$.

Note that $R$ is naturally a subring of $A$ via the map $r \mapsto (r,0)$,
which is a ring homomorphism since the unit of $A$ is $(1,0)$.
Using \zcref{eqn:multA}, one sees that $\epsilon^2 = 0$.
Let us observe that under the assumption that $\sigma$ is not surjective, $\sigma(R)$ is a proper subring of $R$.

Our first goal is to show that under our assumptions, 
the injective dimension of $A$ as a right module over itself is infinite.
To do that, we consider $N = A/(\epsilon\cdot R) \cong R$.
We consider it as a right $A$-module with the action given by
$r \cdot (x+\epsilon y) = r\cdot x$ for $r \in N$, $x+\epsilon y \in A$.

\begin{prop}\label{prop:infinjdim}
Let $R$ be a ring, let $\sigma:R \to R$ be an injective ring homomorphism,
and let $A = R[x;\sigma]/(x^2)$ as in \zcref{eqn:hero}.
Then for any $n \ge 1$,
there is an isomorphism of abelian groups
\[
\opn{Ext}^n_{A^{\opn{op}}}(N,A) \cong R/\sigma(R).
\]
Hence, assuming $\sigma$ is not surjective,
we deduce that $\injdim_A(A_A) = +\infty$.
\end{prop}
\begin{proof}
We will compute these $\opn{Ext}$'s by using a projective resolution of the right $A$-module $N$.
To do that, 
we consider the map $\varphi:A \to A$ given by $\varphi(x+\epsilon y) = \epsilon \cdot (x+\epsilon y) = \epsilon x$. This is a map of right $A$-modules, with the property that both its image and its kernel are equal to $\epsilon R$. Hence,
it gives rise to a projective resolution of the right $A$-module $N$:
\[
\dots \xrightarrow{\varphi} A \xrightarrow{\varphi} A \xrightarrow{\varphi} A \to N \to 0
\]
Using the identification $\opn{Hom}_{A^{\opn{op}}}(A,A) \cong A$,
the map $\opn{Hom}_{A^{\opn{op}}}(\varphi,A)$ is equal to the map $\psi:A \to A$
given by $\psi(x+\epsilon y)= (x+\epsilon y)\cdot \epsilon = \epsilon \sigma(x)$.
Using this projective resolution and this identification, we
thus see that for all $n \ge 1$ it holds that
\[
\opn{Ext}^n_{A^{\opn{op}}}(N,A) \cong \ker(\psi)/\opn{Im}(\psi).
\]
But since $\sigma$ is injective
\[
\ker(\psi) = \{x+\epsilon y \in A \mid \epsilon \sigma(x) = 0\} = \{x+\epsilon y \in A \mid x = 0\} \cong R,
\]
while
\[
\opn{Im}(\psi) = \{\epsilon \sigma(x) \mid x + \epsilon y \in A\} \cong \sigma(R),
\]
which shows that 
\[
\opn{Ext}^n_{A^{\opn{op}}}(N,A) \cong R/\sigma(R).
\]
Finally, assuming that $\sigma$ is not surjective,
we have that $R/\sigma(R) \ne 0$,
so $\opn{Ext}^n_{A^{\opn{op}}}(N,A) \ne 0$ for all $n\ge 1$,
which implies that $\injdim_A(A_A) = +\infty$.
\end{proof}

The next and final goal of this section, which is more difficult, is to show that under certain additional assumptions,
we will have $\injdim_A( \phantom{ }_A A) = 1$.

Following \cite[Definition 2.6]{Nic}, the right conductor of the inclusion of rings $\sigma(R) \subseteq R$ is the left ideal of $R$ given by
\[
(R:\sigma(R))_r = \{x \in R \mid R\cdot x \subseteq \sigma(R)\},
\]
which is, by definition, the largest left ideal of $R$ contained in $\sigma(R)$.

In what follows, we will often have to assume that $1 \in (R:\sigma(R))_r\cdot R$, the right ideal generated by the right conductor.
Note that if $R$ is commutative, 
$(R:\sigma(R))_r\cdot R = (R:\sigma(R))_r$, and $1$ is contained in it if and only if $\sigma$ is surjective, which is impossible in our setup.

More explicitly, the assumption that 
$1 \in (R:\sigma(R))_r\cdot R$ says that there exist elements $r_1,\dots,r_n, s_1,\dots,s_n \in R$, with the property that $R\cdot r_i \subseteq \sigma(R)$ for $1 \le i \le n$,
such that 
\begin{equation}\label{eqn:risi}
\sum_{i=1}^n r_i \cdot s_i = 1.
\end{equation}

\begin{lem}\label{lem:fgProj}
Let $R$ be an associative ring,
let $\sigma:R \to R$ be a ring homomorphism,
and assume that $\sigma$ is injective.
Let $r_1,\dots,r_n, s_1,\dots,s_n \in R$ be elements with the property that
$R\cdot r_i \subseteq \sigma(R)$ for $1 \le i \le n$,
such that $\sum_{i=1}^n r_i \cdot s_i = 1$.
Then the following holds:
\begin{enumerate}
\item As a left module, $R$ is a finitely generated projective
left $\sigma(R)$-module, generated by $s_1,\dots,s_n$.
\item Any homomorphism of left $\sigma(R)$-modules $\varphi:R\to R$
satisfies $\varphi(x) = x\cdot \varphi(1)$ for all $x \in R$.
\end{enumerate}
\end{lem}
\begin{proof}\quad
\begin{enumerate}
\item For any $1 \le i \le n$,
let us set $\phi_i:R \to \sigma(R)$,
$\phi_i(x) = x\cdot r_i$.
Our assumption on the $r_i$ says that indeed this is a map $R \to \sigma(R)$, and clearly it is a map of left $\sigma(R)$-modules.
The assumption $\sum_{i=1}^n r_i \cdot s_i = 1$ implies that
for all $x \in R$ it holds that:
\[
x = x\cdot 1 = x\cdot \sum_{i=1}^n r_i \cdot s_i = 
\sum_{i=1}^n \phi_i(x)\cdot s_i.
\]
By the dual basis lemma, \cite[Lemma 3.5.2]{NCBook},
we deduce that $R$ is projective as a left module over $\sigma(R)$,
and that $s_1,\dots,s_n$ generate it.
\item Let $\varphi:R\to R$ be a homomorphism of left $\sigma(R)$-modules,
and let $x \in R$.
Then
\[
\varphi(x) = \varphi(x\cdot 1) = \varphi\left(x\cdot \sum_{i=1}^n r_i \cdot s_i\right) = \sum_{i=1}^n \varphi((x\cdot r_i) \cdot s_i).
\]
Since for each $1 \le i \le n$ it holds that $x\cdot r_i \in \sigma(R)$,
the left $\sigma(R)$-linearity of $\varphi$ shows that 
$\varphi((x\cdot r_i) \cdot s_i) = (x\cdot r_i) \cdot \varphi(s_i)$,
which shows that
\[
\varphi(x) = x \cdot \sum_{i=1}^n r_i \cdot \varphi(s_i).
\]
In particular, for $x=1$,
we obtain that $\varphi(1) = \sum_{i=1}^n r_i \cdot \varphi(s_i)$,
which shows that $\varphi(x) = x \cdot \varphi(1)$.
\end{enumerate}    
\end{proof}

We now use this lemma to show that this conductor condition allows us to bound the injective dimension of $A$ from the left.

\begin{lem}\label{lem:homIsom}
Let $R$ be a ring, let $\sigma:R \to R$ be an injective ring homomorphism which is not surjective,
and let $A = R[x;\sigma]/(x^2)$ as in \zcref{eqn:hero}.
Suppose that 
$1 \in (R:\sigma(R))_r\cdot R$.
Then there is an isomorphism of left $A$-modules
\[
\opn{Hom}_R\left({}_R A_A, {}_{\sigma} R\right) \cong {}_A A.
\]
\end{lem}
\begin{proof}
Recall that we write elements of $A$ as sums $x+\epsilon y$,
where $x,y \in R$. 
In this notation, let us define a map of left $R$-modules $\varphi:A \to {}_\sigma R$ by $\varphi(x+\epsilon y) = y$.
Indeed, $\varphi$ is left $R$-linear,
since for $r \in R$ we have
\[
\varphi(r\cdot (x+\epsilon y)) = \varphi(r\cdot x + r \cdot \epsilon \cdot y) = \varphi(r\cdot x + \epsilon \sigma(r) \cdot y) = \sigma(r) \cdot y = \sigma(r) \cdot \varphi(x+\epsilon y),
\]
as multiplication by $r$ in ${}_\sigma R$ twists via $\sigma$.

Using $\varphi$, we will define an explicit isomorphism $\Psi: {}_A A \to \opn{Hom}_R\left({}_R A_A, {}_{\sigma} R\right)$ by
\[
\Psi(a)(b) = \varphi(b\cdot a).
\]
First, let us verify that $b \mapsto \varphi(b \cdot a)$ is a map of left $R$-modules. Indeed, for $r \in R$,
we have that $r\cdot b$ maps to
\[
\varphi(r\cdot b \cdot a) = r\cdot \varphi(b\cdot a)
\]
since we already saw that $\varphi$ is $R$-linear.
Next, let us see that $\Psi$ itself is a map of left $A$-modules.
Recall that the left $A$-module structure on 
$\opn{Hom}_R\left({}_R A_A, {}_{\sigma} R\right)$ is given by
$(a\cdot \psi) (b) = \psi(b \cdot a)$ for $a, b \in A$, $\psi \in \opn{Hom}_R\left({}_R A_A, {}_{\sigma} R\right)$.
Then, for $a_1, a_2, b\in A$, 
we have
\[
\Psi(a_1 \cdot a_2)(b) = \varphi(b \cdot a_1 \cdot a_2) = 
\Psi(a_2)(b\cdot a_1) = (a_1 \cdot \Psi(a_2))(b).
\]
Let us show $\Psi$ is injective.
Take $a \in A$ such that $\Psi(a) = 0$,
and write $a = x+\epsilon y$.
The assumption that $\Psi(a) = 0$ implies that $\Psi(a)(b) = 0$ for all $b \in A$.
Taking $b=1$ gives
\[
0 = \Psi(a)(1) = \varphi(a) = y,
\]
while taking $b = \epsilon$ gives
\[
0 = \Psi(a)(\epsilon) = \varphi(\epsilon (x+\epsilon y)) = \varphi(\epsilon \cdot x) = x,
\]
so $a = 0$.
Finally, let us show $\Psi$ is surjective.
Let $\beta \in \opn{Hom}_R\left({}_R A_A, {}_{\sigma} R\right)$,
and let us take $a = \beta(\epsilon) + \epsilon \cdot \beta(1)$.
We will show that $\Psi(a) = \beta$, which will complete the proof.
Let $\alpha = \beta - \Psi(a)$.
Then 
\[
\alpha(1) = \beta(1) - \Psi(a)(1) = \beta(1) - \phi(\beta(\epsilon) + \epsilon \cdot \beta(1)) = \beta(1) - \beta(1) = 0,
\]
and 
\[
\alpha(\epsilon) =  \beta(\epsilon) - \Psi(a)(\epsilon) = 
\beta(\epsilon) - \phi(\epsilon \cdot a) = 
\beta(\epsilon) - \phi(\epsilon \cdot \beta(\epsilon)) = \beta(\epsilon) - \beta(\epsilon) = 0.
\]
Now, taking an arbitrary element $x+\epsilon y \in A$ with $x,y \in R$,
it is enough to show that $\alpha(x) = 0$ and that $\alpha(\epsilon y) =0$
First,
\[
\alpha(x) = \beta(x) - \Psi(a)(x) = \sigma(x)\cdot \beta(1) - \sigma(x)\cdot \Psi(a)(1) = \sigma(x)\cdot \alpha(1) = 0,
\]
as both $\beta$ and $\Psi(a)$ are left $R$-linear where the action on the target is via the $\sigma$ twist.
Define a map $\gamma:R \to R$ by $\gamma(y) = \alpha(\epsilon \cdot y)$. Viewing $R$ as a left $\sigma(R)$-module,
we claim that $\gamma$ is $\sigma(R)$-linear,
because for any $r \in R$,
\[
\gamma(\sigma(r) \cdot y) = \alpha(\epsilon \cdot \sigma(r) \cdot y) = 
\alpha(r\cdot \epsilon \cdot y) = \sigma(r) \cdot \alpha(\epsilon \cdot y) = \sigma(r) \cdot \gamma(y).
\]
According to \zcref{lem:fgProj}(2),
we deduce that for any $y \in R$, 
it holds that
\[
\gamma(y) = y \cdot \gamma(1) = y\cdot \alpha(\epsilon) = 0,
\]
which completes the proof.
\end{proof}

Using the above, we can finally bound the injective dimension of $A$ as a left module over itself:

\begin{prop}\label{prop:injdimleft}
Let $R$ be a ring, let $\sigma:R \to R$ be an injective ring homomorphism which is not surjective,
and let $A = R[x;\sigma]/(x^2)$ as in \zcref{eqn:hero}.
Suppose that 
$1 \in (R:\sigma(R))_r\cdot R$.
Then $\injdim_A({}_A A) \le \injdim_R({}_R R)$.
\end{prop}
\begin{proof}
Note that $\sigma:R \to \sigma(R)$ is a ring isomorphism,
so the forgetful functor $\sigma^*:{}_{\sigma(R)}\opn{Mod} \to {}_R\opn{Mod}$ is an equivalence of categories. 
By \zcref{lem:fgProj}(1), $R$ is a finitely generated projective left $\sigma(R)$-module,
and as being finitely generated projective is preserved by equivalences of module categories, we see that $\sigma^*(R) = {}_\sigma R$ is a finitely generated projective left $R$-module.
Next, observe that as a left $R$-module,
we have that $A \cong R \oplus {}_\sigma R$,
and since both of these are finitely generated projective left $R$-modules,
we deduce that $A$ is also a finitely generated projective left $R$-module.
Let $ n = \injdim_R ({}_R R)$.
We may assume that $n<\infty$, 
for otherwise there is nothing to prove.
Note that since ${}_{\sigma} R$ is a finitely generated projective left $R$-module,
it follows that it is a direct summand of $R^k$ for some $k$,
which implies that $\injdim_R({}_\sigma R) \le \injdim_R({}_R R^k) = \injdim_R({}_R R) = n$.
Take an injective resolution of the left $R$-module ${}_{\sigma} R$:
\[
0 \to {}_{\sigma} R \to I^0 \to I^1 \to \dots \to I^n \to 0.
\]
Let us apply to this exact sequence the functor $\opn{Hom}_R({}_R A_A,-) :{}_R\opn{Mod} \to {}_A \opn{Mod}$. Note that since $A$ is projective as a left $R$-module, this functor is exact, 
so we obtain an exact sequence of left $A$-modules:
\[
0 \to \opn{Hom}_R({}_R A_A,{}_{\sigma} R) \to \opn{Hom}_R({}_R A_A,I^0)  \to \dots \to \opn{Hom}_R({}_R A_A,I^n) \to 0.
\]
For each $0 \le i \le n$,
the left $A$-module $\opn{Hom}_R({}_R A_A,I^i)$ is injective,
so this is an injective resolution of $\opn{Hom}_R({}_R A_A,{}_{\sigma} R)$ of length at most $n$.
By \zcref{lem:homIsom},
we have an isomorphism of left $A$-modules
\[
 \opn{Hom}_R({}_R A_A,{}_{\sigma} R) \cong {}_A A,
\]
which shows that $\injdim_A({}_A A) \le n = \injdim_R({}_R R)$.
\end{proof}

The above is already enough to show that Gorenstein symmetry fails, 
but we wish to be precise about the injective dimension, 
so we go further. Recall that a not necessarily commutative nonzero ring $A$ is called a domain if it has no zero divisors.

\begin{lem}\label{lem:not-inj}
Let $R$ be a domain which is not a division ring, let $\sigma:R \to R$ be an injective ring homomorphism which is not surjective,
and let $A = R[x;\sigma]/(x^2)$ as in \zcref{eqn:hero}.
Then $\injdim_A({}_A A) > 0$.
\end{lem}
\begin{proof}
Since as a right $R$-module $A$ is free of rank $2$,
the functor $A\otimes_R -$ is exact.
If $J$ is any injective left $A$-module,
the isomorphism
\[
\opn{Hom}_R(-,J) \cong \opn{Hom}_A(A\otimes_R -,J)
\]
shows that $J$ is also injective when considered as a left $R$-module.
Assume by contradiction that $A$ is injective as a left $A$-module,
so it is also injective as a left $R$-module.
But $R$, as a left module over itself, is a direct summand of $A$,
so we deduce that $R$ is also self-injective from the left.
But this is impossible for a domain which is not a division ring,
because by Baer's criterion, if $R$ is self-injective,
for any $0\ne a \in R$, the map $R\cdot a\to R$, $r\cdot a \mapsto r$ (which is well defined since $R$ is a domain)
can be extended to a map $R\to R$, which shows $a$ is right-invertible.
Hence, there is $c \in R$ such that $a\cdot c =1$, so that $(c\cdot a-1)\cdot c = 0$ implies $c\cdot a =1$,
contradicting the fact that $R$ is not a division ring.
Hence, $\injdim_A({}_A A) > 0$.
\end{proof}

We finish this section by showing that under our assumptions, if $R$ is two-sided noetherian then so is $A$.

\begin{prop}\label{prop:noeth}
Let $R$ be a ring, let $\sigma:R \to R$ be an injective ring homomorphism which is not surjective,
and let $A = R[x;\sigma]/(x^2)$ as in \zcref{eqn:hero}.
Suppose that 
$1 \in (R:\sigma(R))_r\cdot R$.
Assume that $R$ is both left and right noetherian.
Then $A$ is also left and right noetherian.
\end{prop}
\begin{proof}
We know that $A$ is finitely generated as a right $R$-module by the definition of $A$.
We also know that $A$ is finitely generated as a left $R$-module because both $R$ and ${}_\sigma R$ are finitely generated (the latter by \zcref{lem:fgProj} and the argument in the proof of \zcref{prop:injdimleft}), and as left $R$-modules, $A \cong R \oplus {}_{\sigma} R$.
A ring which is module finite over a right noetherian ring is right noetherian by \cite[Lemma 1.1.3]{NCBook}, which proves the result.
\end{proof}

Recall that a ring $R$ is called left hereditary if its left global dimension is at most $1$.
Here is the main result of this section.

\begin{thm}\label{thm:main-gen-const}
Let $R$ be a two-sided noetherian domain which is left hereditary, and is not a division ring. 
Let $\sigma:R \to R$ be an injective ring homomorphism which is not surjective,
and let $A = R[x;\sigma]/(x^2)$ as in \zcref{eqn:hero}.
Suppose that 
$1 \in (R:\sigma(R))_r\cdot R$.
Then $A$ is a two-sided noetherian ring which satisfies
\[
\injdim_A({}_A A) = 1, \quad \injdim_A(A_A) = +\infty.
\]
\end{thm}
\begin{proof}
By \zcref{prop:noeth}, the ring $A$ is two-sided noetherian.
On the right, according to \zcref{prop:infinjdim} we have that $\injdim_A(A_A) = +\infty$.
On the left, by \zcref{prop:injdimleft}, \zcref{lem:not-inj} and the definition of being hereditary, we obtain:
\[
0 < \injdim_A({}_A A) \le \injdim_R({}_R R) \le 1,
\]
and as these are integers, 
we deduce that $\injdim_A({}_A A) = 1$.
\end{proof}

The next few sections will be dedicated to the construction of a ring $R$ which satisfies the conditions of \zcref{thm:main-gen-const}.

\section{A quaternionic division ring}\label{sec:div}

The purpose of this section is to construct a quaternion-like division ring $D$ with some special properties, which will be used in the next section. It is a variant of a division ring constructed in \cite[Page 13, after Lemma 4.3]{SW}.
There, it is used in \cite[Theorem 4.11]{SW}, 
to construct simple hereditary noetherian domains with prescribed Grothendieck group.

The basic idea of the construction is to take a union of an infinite set of smaller division rings,
each of which is obtained by localizing a quantum torus.

Before constructing it, we state a theorem,
the main result of this section, 
which is essentially a list of properties which we will require from $D$.
These will be the only properties of $D$ used in the remainder of this paper.

\begin{thm}\label{thm:divisionRing}
Let $\k$ be a field of characteristic not $2$.
There exists a division ring $D$ with the following properties:
\begin{enumerate}
\item $D$ is a $\k$-algebra.
\item For any $n \in \mathbb{Z}$,
there are non-zero quaternion-like elements $u_n, v_n \in D$,
such that $u_n\cdot v_n = -v_n\cdot u_n$,
and such that each $u_n^2$ is in the center of $D$.
\item For any finite subset $F \subseteq \mathbb{Z}$,
there is a sub $\k$-algebra $D_F \subseteq D$,
and a subfield $Z_F \subseteq D_F$,
such that $Z_F$ is contained in the center of $D$, and the following hold:
\begin{enumerate}
\item For all $n \in F$, we have that $u_n, u_n^{-1}, v_n, v_n^{-1} \in D_F$.
\item As a vector space, $D_F$ is finite dimensional over $Z_F$.
\item For any $x \in D_F$, and any $n \notin F$, $x\cdot u_n = u_n \cdot x$ and $x\cdot v_n = v_n \cdot x$.
\item If $F_1 \subseteq F_2 \subseteq \mathbb{Z}$ are finite sets of integers,
then $D_{F_1} \subseteq D_{F_2}$.
\item There is an equality $D = \bigcup_F D_F$.
\end{enumerate}
\item There is a $\k$-linear automorphism $\beta:D \to D$ which satisfies $\beta(u_n) = u_{n-1}$ and $\beta(v_n) = v_{n-1}$ for all $n \in \mathbb{Z}$.
\item There is a bijective $\k$-linear function $\iota:D \to D$,
such that $\iota(x\cdot y) = \iota(y)\cdot \iota(x)$ for all $x,y \in D$,
and such that $\iota\circ \iota = 1_D$, and $\iota(u_n) = -u_n$, $\iota(v_n) = v_n$ for all $n \in \mathbb{Z}$.
\end{enumerate} 
\end{thm}

Fix any field $\k$ with characteristic not equal to $2$.
We first build a quantum torus.
In algebraic geometry, the $n$-dimensional torus $(\k^\times)^n$ has  coordinate ring being 
\[
\k[x_1,\dots,x_n,x_1^{-1},\dots,x_n^{-1}].
\]
Now, take a finite subset $F \subseteq \mathbb{Z}$,
which will be the index set for two families of variables.
Following \cite[Example 1.5.10(ii)]{NCBook},
the quantum torus $\mathcal{O}^{2F}_{(\pm 1)}(\k^{\times})$ with coefficients being equal to $\pm 1$ will be defined as the $\k$-algebra generated by elements $u_n ,v_n, u^{-1}_n, v^{-1}_n$ for all $n\in F$,
subject only to the following relations 
\[
u_n \cdot v_m = v_m \cdot u_n, \quad u_n \cdot u_m = u_m \cdot u_n, \quad v_n \cdot v_m = v_m \cdot v_n
\]
for $n,m \in F$ with $n \ne m$, and
\[
u_n \cdot u_n^{-1} = u_n^{-1} \cdot u_n = v_n\cdot v_n^{-1} = v_n^{-1} \cdot v_n = 1
\]
for all $n \in F$,
and most importantly, the quantum relation
\begin{equation}\label{eqn:quantum}
u_n \cdot v_n = -v_n\cdot u_n,
\end{equation}
for all $n \in F$.

The rings $D_F$, whose union will be the division ring $D$ we are trying to construct, will be obtained from the quantum torus by inverting more elements in it. To do that, we identify a large central subring of the quantum torus, which will allow us to localize the quantum torus.

Squaring the quantum torus gives an ordinary torus, in the following sense: let us define $\mathcal{P}^{2F}_{(\pm 1)}(\k^{\times})$ to be the sub $\k$-algebra of $\mathcal{O}^{2F}_{(\pm 1)}(\k^{\times})$ generated by the elements of the form $u_n^{\pm 2}, v_n^{\pm 2}$ for all $n \in F$.
Here are some basic properties of $\mathcal{P}^{2F}_{(\pm 1)}(\k^{\times})$.

\begin{prop}\label{prop:square-torus}
The following hold for $\mathcal{P}^{2F}_{(\pm 1)}(\k^{\times})$:
\begin{enumerate}
\item $\mathcal{P}^{2F}_{(\pm 1)}(\k^{\times})$ is contained in the center of $\mathcal{O}^{2F}_{(\pm 1)}(\k^{\times})$.
\item $\mathcal{O}^{2F}_{(\pm 1)}(\k^{\times})$ is a finitely generated  $\mathcal{P}^{2F}_{(\pm 1)}(\k^{\times})$-module.
\end{enumerate}
\end{prop}
\begin{proof}\quad
\begin{enumerate}
\item It is enough to show that elements of the form $u_n^2$ and $v_n^2$ for $n \in F$ are in the center of $\mathcal{O}^{2F}_{(\pm 1)}(\k^{\times})$ (as their inverses will then automatically also belong to the center). Doing a reduction to generators and inverses again,
it is enough to show these elements commute with each $v_m, u_m$ for all $m \in F$. Each $u_n$ (respectively $v_n$) commutes with $u_m$, for all $m\in F$ (resp. $v_m$, for all $m \in F$).
If $m\ne n$,
then $u_n$ itself already commutes with $v_m$,
so clearly $u_n^2$ commutes with each of $u_m, v_m$ and $v_n^2$ commutes with each of $u_m, v_m$.
Thus, we only need to check that $u_n^2$ (respectively $v_n^2$) commutes with $v_n$ (resp. $u_n$),
and this follows immediately from \zcref{eqn:quantum}.
\item Using the relations between the $u_n$'s and $v_n$'s,
any element of $\mathcal{O}^{2F}_{(\pm 1)}(\k^{\times})$ can be written as a $\k$-linear combination of monomials of the form
\[
\prod_{n \in F} u_n^{a_n} \cdot v_n^{b_n},
\]
where each $a_n, b_n \in \mathbb{Z}$.
Since any integer is either even or odd,
we may write $a_n = 2\cdot \alpha_n + \epsilon_n$ and $b_n = 2\cdot \beta_n + \delta_n$, 
where $\alpha_n, \beta_n \in \mathbb{Z}$, and $\delta_n, \epsilon_n \in \{0,1\}$.
We thus see that we only need to generate the monomials
\[
\prod_{n \in F} u_n^{2\cdot \alpha_n} \cdot v_n^{2\cdot \beta_n} \cdot u_n^{\epsilon_n}\cdot v_n^{\delta_n}.
\]
But as $\prod_{n \in F} u_n^{2\cdot \alpha_n} \cdot v_n^{2\cdot \beta_n} \in \mathcal{P}^{2F}_{(\pm 1)}(\k^{\times})$,
we see that 
\[
\prod_{n \in F} u_n^{\epsilon_n}\cdot v_n^{\delta_n}
\]
for all $\epsilon_n, \delta_n \in \{0,1\}$
generate $\mathcal{O}^{2F}_{(\pm 1)}(\k^{\times})$ as an  $\mathcal{P}^{2F}_{(\pm 1)}(\k^{\times})$-module,
and since $F$ is a finite set,
this is a finite set of generators,
whose size is at most $4^{|F|}$.
\end{enumerate}    
\end{proof}

\begin{rem}
One can show that $\mathcal{P}^{2F}_{(\pm 1)}(\k^{\times})$ is naturally isomorphic, as a $\k$-algebra, to the coordinate ring of the ordinary $2\cdot|F|$-dimensional algebraic torus, $\k[x_i^{\pm 1},y_i^{\pm 1}]_{i \in F}$.
An isomorphism 
\[
\k[x_i^{\pm 1},y_i^{\pm 1}]_{i \in F} \to \mathcal{P}^{2F}_{(\pm 1)}(\k^{\times})
\]
is given 
by $x_i^{\pm 1} \mapsto u_i^{\pm 2}$,
and $y_i^{\pm 1} \mapsto v_i^{\pm 2}$.
This map is clearly surjective,
and its injectivity can be deduced from the fact that the monomials in the $u_i, v_i$ form a $k$-basis of $\mathcal{O}^{2F}_{(\pm 1)}(k^{\times})$, see \cite[Example 1.5.10(ii)]{NCBook}.
\end{rem}

\begin{prop}\label{prop:is-domain}
The ring $\mathcal{O}^{2F}_{(\pm 1)}(\k^{\times})$ is a domain,
and the ring $\mathcal{P}^{2F}_{(\pm 1)}(\k^{\times})$ is an integral domain.
\end{prop}
\begin{proof}
Since by \zcref{prop:square-torus} the ring $\mathcal{P}^{2F}_{(\pm 1)}(\k^{\times})$ is a commutative subring of $\mathcal{O}^{2F}_{(\pm 1)}(\k^{\times})$,
the second claim follows from the first one.
The first claim follows from \cite[Example 1.5.10(ii)]{NCBook}, 
where $\mathcal{O}^{2F}_{(\pm 1)}(\k^{\times})$ is built as an iterated  skew Laurent polynomial ring, and using \cite[Theorem 1.2.9(i), Theorem 1.4.5]{NCBook} repeatedly applied for each variable.
\end{proof}

To shorten notation, in what follows we will further write
$\mathcal{O}_F := \mathcal{O}^{2F}_{(\pm 1)}(\k^{\times})$
and $\mathcal{P}_F := \mathcal{P}^{2F}_{(\pm 1)}(\k^{\times})$.
We denote by $Z_F := \opn{Frac}(\mathcal{P}_F)$ the field of fractions of the integral domain $\mathcal{P}_F$.
Then, we define
\begin{equation}
D_F := \mathcal{O}_F \otimes_{\mathcal{P}_F} Z_F
\end{equation}
for any finite subset $F \subseteq \mathbb{Z}$.
Thus, $D_F$ is obtained from the quantum torus $\mathcal{O}_F$ by inverting all nonzero elements of the squared quantum torus $\mathcal{P}_F$.

The next result summarizes some basic properties of $D_F$:

\begin{prop}\label{prop:DFProp}
The ring $D_F$ satisfies the following properties:
\begin{enumerate}
\item It is a $\k$-algebra.
\item The natural map $\mathcal{O}_F \to D_F$ is injective.
\item The ring $D_F$ is a domain.
\item The ring $Z_F$ is contained in the center of $D_F$.
\item $D_F$ is a finite dimensional vector space over $Z_F$.
\item The ring $D_F$ is a division ring.
\end{enumerate}
\end{prop}
\begin{proof}\quad
\begin{enumerate}
\item This is clear from the definition of $D_F$,
since $\mathcal{O}_F$ is a $\k$-algebra.
\item The map $\mathcal{O}_F \to D_F$ is a localization map,
so its kernel is $\{x \in \mathcal{O}_F \mid \exists 0\ne z \in \mathcal{P}_F, x\cdot z =0\}$, but by \zcref{prop:is-domain}, $\mathcal{O}_F$ is a domain, so this kernel is $0$.
\item By definition,
an element of $D_F$ has the form $x\cdot z^{-1}$,
where $x \in \mathcal{O}_F$, and $0 \ne z \in \mathcal{P}_F$.
Thus, taking two such elements $x_1\cdot z_1^{-1}, x_2 \cdot z_2^{-1}$,
with $x_1,x_2 \in \mathcal{O}_F$ and $0 \ne z_1,z_2 \in \mathcal{P}_F$,
since by \zcref{prop:square-torus} we know that $z_1, z_2$ are in the center of $\mathcal{O}_F$,
we obtain:
\[
x_1\cdot z_1^{-1} \cdot x_2\cdot z_2^{-1} = (x_1 \cdot x_2) \cdot (z_1 \cdot z_2)^{-1}.
\]
Since $\mathcal{O}_F$ is a domain,
this being zero in the localization $D_F$ implies $x_1\cdot x_2 = 0$ in $\mathcal{O}_F$, so that either $x_1\cdot z_1^{-1} = 0$,
or $x_2\cdot z_2^{-1} = 0$, showing $D_F$ is a domain.
\item Take an element $x\cdot z^{-1} \in D_F$, where $x \in \mathcal{O}_F$, and $0 \ne z \in \mathcal{P}_F$.
Next, take $a\cdot b^{-1} \in Z_F$,
with $a,b \in \mathcal{P}_F$ and $b \ne 0$.
Then as $a,b$ and $z$ are in the center of $\mathcal{O}_F$,
the inverses of $b,z$ commute in $D_F$ with elements of $\mathcal{O}_F$,
which implies that
\[
(x\cdot z^{-1}) \cdot (a\cdot b^{-1}) =  (a\cdot b^{-1}) \cdot (x\cdot z^{-1}).
\]
\item This follows from \zcref{prop:square-torus} and the standard commutative algebra fact that if $B$ is a commutative ring, $M$ is a finitely generated $B$-module,
and $S \subseteq B$ is a multiplicatively closed set,
then $S^{-1}M$ is a finitely generated $S^{-1}B$-module.
\item Take a non-zero element $a \in D_F$.
By (3), the map $a\cdot -:D_F \to D_F$ is injective,
and by (4), it is $Z_F$-linear.
Hence, it is an injective map of finite dimensional vector spaces (by (5)),
so it is surjective, which shows that $a$ has a right inverse.
The same argument shows that $-\cdot a:D_F \to D_F$ is surjective,
so $a$ is invertible, and $D_F$ is a division ring.
\end{enumerate}    
\end{proof}

For different finite sets $F \subseteq \mathbb{Z}$,
the rings $D_F$ are related to each other as follows:

\begin{prop}\label{prop:dirLimitconst}
\quad
\begin{enumerate}
\item If $F_1 \subseteq F_2$ are finite subsets of $\mathbb{Z}$,
then there is a unique $\k$-algebra map $\varphi_{1,2}:D_{F_1} \to D_{F_2}$, with $\varphi_{1,2}(u_n) = u_n$ and $\varphi_{1,2}(v_n) = v_n$ for all $n \in F_1$. 
The map $\varphi_{1,2}$ is injective,
and when $F_2 = F_1$, the map $\varphi_{1,1}$ is the identity of $D_{F_1}$.
\item If $F_1 \subseteq F_2 \subseteq F_3$ are two such containments, then $\varphi_{1,3} = \varphi_{2,3}\circ \varphi_{1,2}$.
\end{enumerate}
\end{prop}
\begin{proof}
\quad
\begin{enumerate}
\item We define a $\k$-algebra map $\psi_{1,2}:\mathcal{O}_{F_1} \to \mathcal{O}_{F_2}$ by letting $\psi_{1,2}(u_n^{\pm 1}) = u_n^{\pm 1}$ and $\psi_{1,2}(v_n^{\pm 1}) = v_n^{\pm 1}$ for all $n \in F_1$. 
This defines a unique such map since these elements generate $\mathcal{O}_{F_1}$, 
and since any defining relation between these generators in $\mathcal{O}_{F_1}$ is, by definition, present also in $\mathcal{O}_{F_2}$.
Moreover, $\psi_{1,2}$ is an inclusion of quantum tori,
\cite[1.5.9 and Example 1.5.10(ii)]{NCBook}, so the map $\psi_{1,2}$ is in fact injective.
By definition, the map $\psi_{1,2}$ maps $\mathcal{P}_{F_1}$ into $\mathcal{P}_{F_2}$, and as it is injective,
non-zero such elements are mapped to non-zero elements of $\mathcal{P}_{F_2}$.
By the universal property of localization, it extends uniquely to a map $\varphi_{1,2}$ on the localizations $D_{F_1} \to D_{F_2}$. Since $D_{F_1}$ is a division ring and this map is non-zero, it is injective.
In particular, uniqueness shows that $\varphi_{1,1}$ is the identity.
\item This follows from the corresponding property of $\psi_{i,j}$ and the uniqueness of localization maps.
\end{enumerate}
\end{proof}

The above result shows that the $D_F$ form a direct system of $\k$-algebras.
Taking its colimit, we obtain the main construction of this section:

\begin{equation}
D := \varinjlim_{F \subseteq \mathbb{Z}, |F| < \infty} D_F.
\end{equation}

Since the transition maps $\varphi_{1,2}:D_{F_1} \to D_{F_2}$ are all injective, so are the natural maps $D_F \to D$ for each such $F$.
Hence, we will identify $D_F$ as a sub $\k$-algebra of $D$.
This also immediately gives us:

\begin{prop}\label{prop:dProperties}
The algebra $D$ satisfies the following properties:
\begin{enumerate}
\item Any element of $D$ lies in $D_F$ for some finite subset $F \subseteq \mathbb{Z}$.
\item The $\k$-algebra $D$ is a division ring.
\end{enumerate}
\end{prop}
\begin{proof}\quad
\begin{enumerate}
\item This is because the colimit defining $D$ is over a directed set.
\item
If $x \in D$ is any non-zero element,
take some finite subset $F \subseteq \mathbb{Z}$,
such that $0\ne x \in D_F$.
Since $D_F$ is a division ring, $x$ has an inverse $y$,
which also belongs to $D$,
so any non-zero element of $D$ is invertible.
\end{enumerate}
\end{proof}

We need one final property of $D$ before proving the main result,
a universal property that $D$ satisfies:
\begin{prop}\label{prop:univerPropD}
Let $E$ be any division ring which is a $\k$-algebra,
and assume that for any $n \in \mathbb{Z}$,
there are non-zero elements $a_n, b_n \in E$ which satisfy:
\begin{enumerate}
\item For all $n \in \mathbb{Z}$, $a_n \cdot b_n = -b_n \cdot a_n$.
\item For all $n,m \in \mathbb{Z}$ with $n\ne m$,
it holds that
\[
a_n \cdot b_m = b_m \cdot a_n, \quad a_n \cdot a_m = a_m \cdot a_n, \quad b_n \cdot b_m = b_m \cdot b_n.
\]
\end{enumerate}
Then the following holds:
\begin{enumerate}
\item For any finite subset $F \subseteq \mathbb{Z}$,
there is a unique $\k$-algebra map $\alpha_F:\mathcal{O}_F \to E$ which satisfies $\alpha_F(u_n) = a_n$ and $\alpha_F(v_n) = b_n$ for all $n \in F$.
\item Suppose that for any finite subset $F \subseteq \mathbb{Z}$,
the restriction of $\alpha_F$ to the commutative subring $\mathcal{P}_F \subseteq \mathcal{O}_F$ is injective.
Then there is a unique $\k$-algebra map $\alpha:D \to E$ such that $\alpha(u_n) = a_n$ and $\alpha(v_n) = b_n$ for all $n \in \mathbb{Z}$.
\end{enumerate}
\end{prop}
\begin{proof}\quad
\begin{enumerate}
\item This is immediate from the universal property of the presentation of $\mathcal{O}_F$ in terms of generators and relations.
\item First, we show uniqueness.
Assume $\alpha, \beta:D \to E$ are two such maps,
and consider $T := \{x \in D \mid \alpha(x) = \beta(x)\}$.
Then $T$ is a $\k$-subalgebra of $D$ which by definition contains all the $u_n, v_n$ for all $n \in \mathbb{Z}$.
But since a map of division rings preserves inverses,
we deduce that $T$ is also a division ring.
Finally, since $D$ is generated as a division ring by $\k$ and the $u_n$'s and $v_n$'s (because that is true for the $D_F$, and any element of $D$ belongs to some $D_F$), we deduce that $T = D$, so $\alpha = \beta$.
Next, we show existence. 
Since the map $\alpha_F$ is injective on $\mathcal{P}_F$,
and since $E$ is a division ring,
localization allows us to extend it to a map $D_F \to E$,
which we will also denote by $\alpha_F$.
Moreover, using the same argument as in the uniqueness part of the proof,
for $F_1 \subseteq F_2$,
and after identifying each $D_F$ with its image in $D$,
we have that $\alpha_{F_2}$ is an extension of $\alpha_{F_1}$.
Hence, passing to the direct limit along all finite subsets $F \subseteq \mathbb{Z}$,
we obtain the required $\alpha$.
\end{enumerate}    
\end{proof}

We are now ready to prove that the division ring $D$ satisfies the conditions of \zcref{thm:divisionRing},
the main result of this section:
\begin{proof}\quad
\begin{enumerate}
\item This is because $D$ is a direct limit of $\k$-algebras,
and the transition maps are $\k$-algebra maps.
\item As we identified the $D_F$'s as subalgebras of $D$,
we have for each $n \in \mathbb{Z}$ non-zero elements $u_n, v_n \in D$ such that $u_n \cdot v_n = -v_n\cdot u_n$. It remains to show that $u_n^2$ is in the center of $D$.
Let $x \in D$ be an arbitrary element.
By \zcref{prop:dProperties},
there exists a finite subset $F_1 \subseteq \mathbb{Z}$,
such that $x \in D_{F_1}$.
Let $F = F_1 \cup \{n\}$.
Then $x, u_n^2 \in D_F$,
and since $u_n^2 \in Z_F$, by \zcref{prop:DFProp}(4),
we know that $u_n^2$ is in the center of $D_F$.
Hence, $x\cdot u_n^2 = u_n^2\cdot x$,
and this relation remains true in $D$,
showing that $u_n^2$ is in the center of $D$.
\item We have already established the existence of $Z_F \subseteq D_F \subseteq D$. The same argument as in (2) shows that $Z_F$ is contained in the center of $D$. We now prove the rest of the claims in this item:
\begin{enumerate}
\item We have seen that for all $n \in F$, we have that $u_n,u_n^{-1},v_n,v_n^{-1} \in D_F$.
\item The fact that $D_F$ is finite dimensional over $Z_F$ is \zcref{prop:DFProp}(5).
\item If $x \in D_F$ and $n \notin F$,
then we may write $x$ as a finite sum
\[
x = \sum a_i \cdot z_i^{-1} \cdot \prod_{j \in F} u_j^{m_j}  \cdot \prod_{k \in F}  v_k^{l_k}  
\]
where $a_i\in \k$, $m_j, l_k \in \mathbb{Z}$,
and $z_i \in Z_F$.
But as $n \notin F$,
each of the $u_j^{m_j}$ and of the $v_k^{l_k}$ commute with $u_n$.
As for $z_i$, 
we have that $z_i \in Z_F \subseteq Z_{F \cup \{n\}}$,
which is contained in the center of $D_{F\cup\{n\}}$,
so $z_i$ commutes with $u_n$,
and hence $z_i^{-1}$ also commutes with $u_n$.
Hence, we deduce that $x\cdot u_n = u_n \cdot x$.
The same argument also shows that $x\cdot v_n = v_n \cdot x$.
\item The fact that for $F_1 \subseteq F_2$ we have $D_{F_1} \subseteq D_{F_2}$ is \zcref{prop:dirLimitconst}(1) together with the identification of the $D_F$'s with their image inside $D$.
\item The equality $D = \bigcup_F D_F$ is our definition of $D$.
\end{enumerate}
\item We will use the universal property of \zcref{prop:univerPropD} to construct a $\k$-linear automorphism $\beta:D \to D$ which satisfies $\beta(u_n) = u_{n-1}$ and $\beta(v_n) = v_{n-1}$ for all $n \in \mathbb{Z}$.
To do that, apply \zcref{prop:univerPropD} with $E = D$,
$a_n = u_{n-1}$ and $b_n = v_{n-1}$ for all $n \in \mathbb{Z}$.
For a finite subset $F \subseteq \mathbb{Z}$,
we define $\alpha_F :\mathcal{O}_F \to E$ as follows.
First, let $F' = \{n-1 \mid n \in F\}$.
Next, let $\alpha_F':\mathcal{O}_F \to \mathcal{O}_{F'}$ be the map
which satisfies $\alpha_F'(u_n) = u_{n-1}$ and $\alpha_F'(v_n) = v_{n-1}$. Because the generators and relations of $\mathcal{O}_F$ and $\mathcal{O}_{F'}$ are identical up to change of names,
it is clear that $\alpha_F'$ is a $\k$-algebra isomorphism.
Then, let $\alpha_F$ be the composition of $\alpha_F'$ with the inclusion $\mathcal{O}_{F'} \inj E$. Then it is injective,
and in particular its restriction to $\mathcal{P}_F$ is injective.
Moreover, by construction, it satisfies that $\alpha_F(u_n) = a_n$ and $\alpha_F(v_n) = b_n$ for all $n\in F$.
Hence, by \zcref{prop:univerPropD},
we obtain a $\k$-algebra map $\beta:D \to E=D$ which satisfies that $\beta(u_n)=u_{n-1}$ and $\beta(v_n) = v_{n-1}$ for all $n\in \mathbb{Z}$.
Repeating the same argument with $+1$ replacing $-1$ will provide us with a $\k$-algebra map $\gamma:D \to D$ which satisfies that $\gamma(u_n)=u_{n+1}$ and $\gamma(v_n) = v_{n+1}$ for all $n\in \mathbb{Z}$.
Finally, each of $\beta\circ \gamma$ and $\gamma\circ \beta$ is a $\k$-algebra map of division rings which is the identity on the $u_n$ and the $v_n$, so, as in the proof of \zcref{prop:univerPropD},
we deduce that $\beta, \gamma$ are inverses of each other, so are automorphisms.
\item We wish to construct a bijective $\k$-linear function $\iota:D \to D$,
such that $\iota(x\cdot y) = \iota(y)\cdot \iota(x)$ for all $x,y \in D$,
and such that $\iota\circ \iota = 1_D$, and $\iota(u_n) = -u_n$, $\iota(v_n) = v_n$ for all $n \in \mathbb{Z}$.
To do that, we again use the universal property \zcref{prop:univerPropD}.
Let $E = D^{\opn{op}}$,
let $a_n = -u_n$, and let $b_n = v_n$ for all $n \in \mathbb{Z}$.
Observe that this data still satisfies the conditions of the universal property of $D$.
For any finite subset $F \subseteq \mathbb{Z}$,
we define $\alpha_F : \mathcal{O}_F \to E$,
as follows. 
First, consider the map $\mathcal{O}_F \to \mathcal{O}_F^{\opn{op}}$ given by $u_n \mapsto -u_n$ and $v_n \mapsto v_n$. 
This map is clearly bijective. Composing it with the inclusion $\mathcal{O}_F^{\opn{op}} \inj E$ gives an injective map $\alpha_F$.
In particular it is injective on $\mathcal{P}_F$.
Hence, by \zcref{prop:univerPropD},
we obtain an induced $\k$-algebra map $\iota:D \to E$ which satisfies that $\iota(u_n) = -u_n$ and $\iota(v_n) = v_n$ for all $n \in \mathbb{Z}$. 
Considered as a map $\iota:D\to D$,
we deduce that $\iota$ is a $\k$-linear map,
and that $\iota(x\cdot y) = \iota(y)\cdot \iota(x)$ for all $x,y \in D$.
Finally, since $\iota\circ \iota$ is a $\k$-algebra map of division rings which is the identity on the $u_n$ and the $v_n$,
we deduce it is the identity on the whole of $D$.
This implies $\iota$ is bijective.
\end{enumerate}
\end{proof}

\section{An intersection of idealizers}\label{sec:R}

In this section we will construct a ring $R$ which will be shown to satisfy the conditions of \zcref{thm:main-gen-const}.
Throughout this section, fix a field $\k$ of characteristic not $2$,
and let $D$ be a division ring which is a $\k$-algebra that possesses the properties of \zcref{thm:divisionRing}.

\subsection{The ring $R$}
Let $S = D[t]$ be the polynomial algebra over $D$,
with $t$ being the indeterminate, in particular, $t$ is a central element.
A good reference for polynomial rings over division rings is \cite[Section 16]{Lam}. 
Our target ring $R$ will be a certain sub $\k$-algebra of $S$.

Somewhat unusually, and working with the left-right mirror of the convention in \cite{Lam},
given a polynomial $f \in S$,
which is, by definition, a finite sum
\[
f = \sum_{k=0}^n t^k \cdot a_k \in S,
\]
with each $a_k \in D$, 
and given some $b \in D$,
we will consider the left evaluation of $f$ at $b$,
which is defined to be the element
\[
f_{\ell}(b) = \sum_{k=0}^n b^k \cdot a_k \in D.
\]

Here is an elementary fact about this notion we will frequently use:
\begin{lem}\label{lem:belongTo}
Let $f \in S$ and let $b \in D$.
Then the following holds:
\begin{enumerate}
\item It holds that $f-f_{\ell}(b)\in (t-b)\cdot S$.
\item We have that $(t-b)\cdot S \cap D = 0$.
\end{enumerate}    
\end{lem}
\begin{proof}
\quad
\begin{enumerate}
\item This follows immediately from the left version of the first part of \cite[Proposition 16.2]{Lam}.
\item This follows immediately from the left version of the second part of \cite[Proposition 16.2]{Lam}.
\end{enumerate}
\end{proof}

For any $n\in \mathbb{Z}$,
we shall often need the following elements of $S$, 
so we give them names. Let $g_n = t-u_n \in S$ and its conjugate $h_n = t+u_n \in S$.
Denote by $C_D(u_n)$ the centralizer of $u_n$ in $D$.
By definition,
\[
C_D(u_n) = \{x \in D \mid x\cdot u_n = u_n \cdot x\}.
\]
This is a sub $\k$-algebra of $D$ which contains the center of $D$.
For any $n \in \mathbb{Z}$,
let 
\begin{equation}\label{eqn:idealizer}
\mathcal{I}_n := \{f \in S \mid f_{\ell}(u_n) \in C_D(u_n)\},
\end{equation}
or more explicitly
\[
\mathcal{I}_n = \{f \in S \mid f_{\ell}(u_n) \cdot u_n = u_n \cdot f_{\ell}(u_n)\}.
\]
Alternatively, one can define $\mathcal{I}_n$ using the notion of idealizers:
Recall, following \cite[1.1.11 and Section 5.5]{NCBook} that if $B$ is a ring, and $J$ is a right ideal in $B$
the idealizer of $J$ in $B$ is the subring $\mathcal{I}_B(J)$ of $B$ given by
\[
\mathcal{I}_B(J) = \{b \in B \mid b\cdot J \subseteq J\}.
\]
Thus, the idealizer is the largest subring of $B$ in which $J$ becomes a two-sided ideal.
Going back to our context, we have the following result which connects these definitions:

\begin{prop}\label{prop:idealizer}
The following holds for all $n \in \mathbb{Z}$:
\begin{enumerate}
\item The set $\mathcal{I}_n$ defined in \zcref{eqn:idealizer} is equal to the idealizer of the right ideal $g_n\cdot S$ in $S$.
\item The set $\mathcal{I}_n$ is a sub $\k$-algebra of $S$ which contains the center of $S$,
and $g_n\cdot S \subseteq \mathcal{I}_n$ is a two-sided ideal in it.
\end{enumerate}
\end{prop}
\begin{proof}
\quad
\begin{enumerate}
\item Let $J = g_n \cdot S$,
which is a right ideal of $S$.
An element $f \in S$ belongs to the idealizer of $J$, 
$\mathcal{I}_S(J)$
if and only if $f\cdot g_n \cdot S \subseteq g_n \cdot S$,
if and only if 
\[
f\cdot g_n \in g_n \cdot S = (t-u_n)\cdot S.
\]
By \zcref{lem:belongTo},
this is equivalent to $(f\cdot g_n)_{\ell}(u_n) = 0$.
Write $f = \sum_{k=0}^m t^k\cdot a_k$.
Then 
\[
f \cdot g_n = f\cdot (t-u_n) = \sum_{k=0}^m t^{k+1}\cdot a_k - \sum_{k=0}^m t^k\cdot a_k \cdot u_n.
\]
Plugging $u_n$, we obtain that
\[
(f\cdot g_n)_{\ell}(u_n) = 
\sum_{k=0}^m u_n^{k+1}\cdot a_k - \sum_{k=0}^m u_n^k\cdot a_k \cdot u_n = 
u_n \cdot f_{\ell}(u_n) - f_{\ell}(u_n)\cdot u_n.
\]
Thus, we see that $f \in \mathcal{I}_S(J)$ 
if and only if $f_{\ell}(u_n) \in C_D(u_n)$,
which is exactly the set $\mathcal{I}_n$ defined in \zcref{eqn:idealizer}.
\item This follows from (1), since the idealizer satisfies these properties.
\end{enumerate}
\end{proof}

We may now define our target $\k$-algebra $R$,
as an intersection of these idealizers.
Note that while we defined $\mathcal{I}_n$ for all $n \in \mathbb{Z}$,
in the definition below, we only take the intersection of these for $n \ge 1$. We define:
\begin{equation}\label{eqn:theR}
R := \bigcap_{n=1}^{\infty} \mathcal{I}_n.    
\end{equation}

\begin{prop}\label{prop:rbasic}
The object $R$ defined in \zcref{eqn:theR} is a sub $\k$-algebra of $S$.
Moreover, $R$ is a domain which is not a division ring,
and it contains the center of $S$.
\end{prop}
\begin{proof}
Since by \zcref{prop:idealizer}(2), each $\mathcal{I}_n$ is a sub $\k$-algebra of $S$, it follows that $R$ is also a sub $\k$-algebra of $S$.
As $S$ is a polynomial ring over a division ring,
it is a domain, so $R$, being a subring of it,
is also a domain.
By \zcref{prop:idealizer}(2), we know that each $\mathcal{I}_n$ contains the center of $S$,
so $R$ also contains the center of $S$.
Finally, $t \in S$ belongs to the center of $S$ by definition,
so $t \in R$, but $t$ is not even invertible in $S$,
so it is also not invertible in $R$. Thus, $R$ is not a division ring.
\end{proof}

\subsection{The endomorphism $\sigma:R \to R$}\label{sec:sigmaDef}

Our next goal is to construct an injective ring homomorphism $\sigma:R \to R$ which is not surjective, and which satisfies $1 \in (R:\sigma(R))_r\cdot R$ as required by \zcref{thm:main-gen-const}.
We will further ensure that $\sigma$ is a $\k$-algebra homomorphism,
which implies that our final construction of $A$, as in \zcref{sec:gen},  will also be a $\k$-algebra.

In the remainder of this section,
let $\beta:D \to D$ be a $\k$-linear automorphism which satisfies $\beta(u_n) = u_{n-1}$ and $\beta(v_n) = v_{n-1}$ for all $n \in \mathbb{Z}$. Such a $\beta$ exists by \zcref{thm:divisionRing}(4).
We extend $\beta$ to a map $S\to S$,
also denoted by $\beta$,
by defining
\begin{equation}\label{eqn:beta}
\beta\left(\sum_{k=0}^n t^k\cdot a_k\right) = \sum_{k=0}^n t^k \cdot \beta(a_k), \quad a_k \in D.
\end{equation}
Then it is clear that $\beta:S \to S$ is a $\k$-linear automorphism of $S$.

\begin{lem}\label{lem:betaofev}
For any $f \in S$ and any $d\in D$,
the automorphism $\beta:S \to S$ defined in \zcref{eqn:beta} satisfies
\[
\Bigl(\beta(f)_{\ell}\Bigr) (\beta(d)) = \beta(f_{\ell}(d)).
\]
\end{lem}
\begin{proof}
Assume $f = \sum_{k=0}^n t^k\cdot a_k$ with each $a_k \in D$.
Then, we calculate:
\[
\Bigl(\beta(f)_{\ell}\Bigr) (\beta(d)) = 
\sum_{k=0}^n \beta(d)^k\cdot \beta(a_k).
\]
Since $\beta$ is a ring homomorphism,
it commutes with addition and multiplication,
so we have:
\[
\sum_{k=0}^n \beta(d)^k\cdot \beta(a_k) = 
\beta\Bigl(\sum_{k=0}^n d^k \cdot a_k\Bigr) = \beta\Bigl(f_{\ell}(d)\Bigr).
\]
\end{proof}

\begin{lem}\label{lem:betaCent}
For any $n \in \mathbb{Z}$,
the automorphism $\beta:D \to D$ satisfies
\[
\beta\bigl(C_D(u_n)\bigr) = C_D(u_{n-1}).
\]
\end{lem}
\begin{proof}
If $x \in C_D(u_n)$,
then by definition $x\cdot u_n = u_n\cdot x$.
Applying $\beta$ to this equation,
and using the fact that $\beta(u_n) = u_{n-1}$ and that $\beta$ is a ring homomorphism,
we see that $\beta(x) \cdot u_{n-1} = u_{n-1} \cdot \beta(x)$,
so $\beta(x) \in C_D(u_{n-1})$.
Conversely, if $x \in C_D(u_{n-1})$,
as $\beta$ is an isomorphism, taking $y=\beta^{-1}(x)$ and repeating the same argument with $\beta$ replaced by $\beta^{-1}$, we obtain that $y \in C_D(u_n)$.
\end{proof}

We connect the automorphism $\beta$ to the idealizers $\mathcal{I}_n$, defined in \zcref{eqn:idealizer}, in the next result:

\begin{lem}\label{lem:betaIn}
For any $n \in \mathbb{Z}$,
it holds that $\beta(\mathcal{I}_n) = \mathcal{I}_{n-1}$.
\end{lem}
\begin{proof}
By definition,
for any $n \in \mathbb{Z}$,
$f \in \mathcal{I}_n$ if and only if $f_{\ell}(u_n) \in C_D(u_n)$.
By \zcref{lem:betaCent},
this is equivalent to $\beta(f_{\ell}(u_n)) \in C_D(u_{n-1})$.
According to \zcref{lem:betaofev} and the definition of $\beta$, we have
\[
\beta(f_{\ell}(u_n)) = \bigl(\beta(f)\bigr)_{\ell}(\beta(u_n)) = \bigl(\beta(f)\bigr)_{\ell}(u_{n-1}).
\]
Thus, we see that $f \in \mathcal{I}_n$ if and only if
$\bigl(\beta(f)\bigr)_{\ell}(u_{n-1}) \in C_D(u_{n-1})$,
which, by definition, is equivalent to $\beta(f) \in \mathcal{I}_{n-1}$.
Since $\beta$ is bijective, this proves the result.
\end{proof}

From this result we deduce:
\begin{prop}\label{prop:injnotsurj}
The automorphism $\beta:S \to S$ defined in \zcref{eqn:beta} restricts to a $\k$-algebra homomorphism $\beta|_R : R \to R$ which is injective but not surjective.
\end{prop}
\begin{proof}
Let $x \in R$.
By \zcref{eqn:theR}, this is equivalent to $x \in \mathcal{I}_n$ for all $n \ge 1$.
For any $n\ge 1$, $n+1 \ge 1$, so $x \in \mathcal{I}_{n+1}$.
This implies by \zcref{lem:betaIn},
that $\beta(x) \in \mathcal{I}_n$.
Since this holds for all $n \ge 1$,
we deduce that $\beta(x) \in R$.
Since $\beta:S \to S$ was an injective $\k$-algebra homomorphism,
it is still an injective $\k$-algebra homomorphism when restricted to $R$.
It remains to show that $\beta|_R:R \to R$ is not surjective.
To see this, consider
the constant polynomial $f = v_0 + t\cdot 0$.
As it is constant, for any $n$, we have $f_{\ell}(u_n) = v_0$.
And for $n \ne 0$,
by \zcref{thm:divisionRing}(3c),
$v_0$ commutes with $u_n$.
Hence, for all $n \ge 1$, 
we have that $f \in \mathcal{I}_n$, and hence, $f \in R$.
However, since $0 \ne u_0\cdot v_0 = -v_0\cdot u_0$, and as $\opn{char}(\k) \ne 2$, $v_0$ does not commute with $u_0$.
We thus see that $f \notin \mathcal{I}_0$.
As $\beta:S \to S$ is an automorphism,
there is a unique $x \in S$ such that $\beta(x) = f$.
The fact that $f \notin \mathcal{I}_0$ implies by \zcref{lem:betaIn}
that $x \notin \mathcal{I}_1$, and hence $x\notin R$.
Thus, $f\notin \opn{Im}(\beta|_R)$,
so $\beta|_R$ is not surjective.
\end{proof}

We define
\begin{equation}\label{eqn:sigma}
\sigma : R \to R, \quad \sigma:=\beta|_R.
\end{equation}

\subsection{The conductor condition}\label{sec:conductor}
Our next goal is to show that $1 \in (R:\sigma(R))_r\cdot R$.
To do that, we first describe the image of $\sigma$ explicitly:
\begin{prop}\label{prop:imsigma}
The $\k$-algebra homomorphism $\sigma:R \to R$ defined in \zcref{eqn:sigma} satisfies:
\[
\sigma(R) = R \cap \mathcal{I}_0,
\]
where $\mathcal{I}_0$ is the idealizer defined in \zcref{eqn:idealizer}.
\end{prop}
\begin{proof}
If $x \in R$,
then by the definition of $R$, for all $n \ge 1$, $x \in \mathcal{I}_n$,
and hence by \zcref{lem:betaIn}, it holds that $\sigma(x) \in \mathcal{I}_{n-1}$.
Thus, 
\[
\sigma(x) \in \bigcap_{n=0}^{\infty} \mathcal{I}_n = \left(\bigcap_{n=1}^{\infty} \mathcal{I}_n\right) \cap \mathcal{I}_0 = R\cap \mathcal{I}_0.
\]
Conversely, if $y \in R\cap \mathcal{I}_0$,
the above calculation shows that $y \in \mathcal{I}_n$ for all $n \ge 0$,
and letting $x = \beta^{-1}(y) \in S$,
applying \zcref{lem:betaIn} again, 
we see that $x \in \mathcal{I}_n$ for all $n \ge 1$,
so that $x \in R$, and $\sigma(x) = \beta(x) = y$.
\end{proof}

Recall that $h_0 = t + u_0 \in S$.
Next, we define four additional elements of $S$ which will be used in the calculations that will follow:
\begin{equation}\label{eqn:rs-elements}
r_1 = \frac{1}{2}\cdot u_0^{-1}\cdot h_0, \quad s_1 = 1, \quad r_2 = -\frac{1}{2}\cdot u_0^{-1}\cdot v_0^{-1}\cdot h_0, \quad s_2 = v_0.
\end{equation}
These will play the role of the $r_i$'s and the $s_i$'s, as in \zcref{eqn:risi}. First, to show these elements belong to $R$,
we prove a general result about belonging to the idealizers $\mathcal{I}_n$:
\begin{lem}\label{lem:coeffBelong}
Let $F \subseteq \mathbb{Z}$ be a finite set,
let $n \in \mathbb{Z}$, and suppose that $n \notin F$.
Let $f = \sum_{k=0}^m t^k\cdot a_k \in S$,
and suppose that for any $0\le k \le m$, 
it holds that $a_k \in D_F$.
Then $f \in \mathcal{I}_n$.
\end{lem}
\begin{proof}
This is because by \zcref{thm:divisionRing}(3)(c),
since $n\notin F$, elements of $D_F$ commute with $u_n$,
and clearly if $a_k$ commutes with $u_n$, 
then so does $u_n^k\cdot a_k$.
\end{proof}

Taking $F = \{0\}$, we have that the coefficients of $r_1,s_1,r_2,s_2$ belong to $D_F$, so the lemma implies that $r_1,s_1,r_2,s_2 \in \mathcal{I}_n$ for all $n \ne 0$,
and in particular, 
\begin{equation}\label{eqn:risiinR}
r_1,s_1,r_2,s_2 \in R.
\end{equation}

\begin{lem}\label{lem:Sh_0}
The element $h_0 = t + u_0 \in S$ satisfies $S\cdot h_0 \subseteq \mathcal{I}_0$.
\end{lem}
\begin{proof}
Take any $f \in S$.
We may write $f = a_0 + t\cdot a_1 + \dots + t^n \cdot a_n$, with $a_i \in D$.
Let us set $c = f_{\ell}(u_0)$.
To check if $f\cdot h_0 \in \mathcal{I}_0$ we calculate:
\begin{align*}
(f\cdot h_0)_{\ell}(u_0) &= 
\left(\bigl(\sum_{k=0}^n t^k\cdot a_k\bigr)\cdot (t + u_0)\right)_{\ell}(u_0) \\ &= 
\left( \sum_{k=0}^n t^{k+1}\cdot a_k + \sum_{k=0}^n t^k\cdot a_k \cdot u_0 \right)_{\ell}(u_0) \\ &= 
\sum_{k=0}^n u_0^{k+1}\cdot a_k + \sum_{k=0}^n u_0^k\cdot a_k \cdot u_0 = u_0\cdot c + c\cdot u_0.
\end{align*}
Then $f\cdot h_0 \in \mathcal{I}_0$ if and only if $u_0 \cdot c + c\cdot u_0$ commutes with $u_0$, and indeed:
\[
u_0 \cdot (u_0 \cdot c + c\cdot u_0) = u_0^2 \cdot c + u_0\cdot c \cdot u_0,
\]
while 
\[
(u_0 \cdot c + c\cdot u_0) \cdot u_0 = u_0\cdot c \cdot u_0 + c\cdot u_0^2.
\]
Since by \zcref{thm:divisionRing}(2), the element $u_0^2$ is in the center of $D$, the result follows.
\end{proof}

\begin{lem}\label{lem:Rrinsigma}
The elements $r_1, r_2$ defined in \zcref{eqn:rs-elements},
satisfy
\[
R\cdot r_1 \subseteq \sigma(R), \quad R\cdot r_2 \subseteq \sigma(R).
\]
\end{lem}
\begin{proof}
If $r \in \{r_1,r_2\}$,
by their definition, we have that $r \in D\cdot h_0$.
Hence,
\[
R \cdot r \subseteq R \cdot D \cdot h_0 \subseteq S \cdot h_0 \subseteq \mathcal{I}_0, 
\]
where the last containment is by \zcref{lem:Sh_0}.
On the other hand, since $R$ is a ring and by \zcref{eqn:risiinR},
$r \in R$, clearly $R\cdot r \subseteq R$,
so we obtain that
\[
R \cdot r \subseteq R \cap \mathcal{I}_0 = \sigma(R),
\]
where the last equality is \zcref{prop:imsigma}.
\end{proof}

We may now achieve our final goal concerning $\sigma$:
\begin{prop}\label{prop:conductor}
The elements $r_1,s_1,r_2,s_2$ defined in \zcref{eqn:rs-elements} satisfy the following properties:
\begin{enumerate}
\item $r_1,s_1,r_2,s_2 \in R$.
\item $R\cdot r_1 \subseteq \sigma(R)$ and $R\cdot r_2 \subseteq \sigma(R)$.
\item There is an equality $r_1\cdot s_1 + r_2 \cdot s_2 = 1$.
\end{enumerate}
Hence, in the notation of the discussion preceding \zcref{eqn:risi},
it holds that 
\[
1 \in (R:\sigma(R))_r\cdot R.
\]
\end{prop}
\begin{proof}
(1) was established in \zcref{eqn:risiinR},
and (2) was established in \zcref{lem:Rrinsigma}.
As for (3), we calculate:
\begin{align*}
r_1\cdot s_1 + r_2 \cdot s_2 &= 
\frac{1}{2}\cdot u_0^{-1}\cdot h_0 -\frac{1}{2}\cdot u_0^{-1}\cdot v_0^{-1}\cdot h_0 \cdot v_0 \\ &= \frac{1}{2}\cdot u_0^{-1} \cdot \left(h_0 - v_0^{-1}\cdot h_0 \cdot v_0 \right).
\end{align*}
Since $h_0 = t + u_0$,
and as $t$ is in the center of $R$ and $u_0$ and $v_0$ anti-commute,
we have that
\[
v_0^{-1}\cdot h_0 \cdot v_0 = v_0^{-1}\cdot (t + u_0) \cdot v_0 = 
t - u_0,
\]
which shows that 
\[
r_1\cdot s_1 + r_2 \cdot s_2 = 
\frac{1}{2}\cdot u_0^{-1} \cdot \left( t+u_0 -t + u_0\right) = 1,
\]
as claimed.
Finally, the fact that this implies that $1 \in (R:\sigma(R))_r\cdot R$ was shown in the discussion preceding \zcref{eqn:risi}.
\end{proof}

\subsection{Symmetry of the ring $R$}\label{sec:symm}

Our next goal will be to show that $R \cong R^{\opn{op}}$ as $\k$-algebras. This will be a non-trivial isomorphism, not the identity map, as $R$ is of course not commutative. This will later allow us to deduce directly that $R$ is left-noetherian from the fact (which we will prove) that it is right-noetherian.

To construct such an isomorphism,
recall that by \zcref{thm:divisionRing}(5),
there is a bijective $\k$-linear function $\iota:D \to D$,
such that $\iota(x\cdot y) = \iota(y)\cdot \iota(x)$ for all $x,y \in D$,
and such that $\iota\circ \iota = 1_D$, and $\iota(u_n) = -u_n$, $\iota(v_n) = v_n$ for all $n \in \mathbb{Z}$.
As before, we extend $\iota$ to a map $S\to S$, also denoted by $\iota$,
by letting
\begin{equation}\label{eqn:iotaExtended}
\iota\left(\sum_{k=0}^n t^k \cdot a_k\right) = \sum_{k=0}^n t^k \cdot \iota(a_k).
\end{equation}
Then, it is clear that it is still the case that $\iota:S \to S$ is $\k$-linear, and that $\iota\circ \iota = 1_S$, so in particular, $\iota$ is bijective.
To check that $\iota(x\cdot y) = \iota(y)\cdot \iota(x)$ for all $x,y \in S$, it is enough to check on monomials:
if $x = t^n \cdot a$, and $y = t^m\cdot b$, with $a,b \in D$,
then since $t$ is central,
\[
\iota(x\cdot y) = \iota(t^n\cdot a \cdot t^m \cdot b) = 
t^{n+m} \cdot \iota(a\cdot b) = t^{n+m} \cdot \iota(b)\cdot \iota(a) = 
\iota(t^m\cdot b) \cdot \iota(t^n \cdot a) = \iota(y) \cdot \iota(x).
\]

Recall that $g_n = t-u_n$ and $h_n = t+u_n$.
Calculating $\iota$ on these elements gives:
\begin{equation}\label{eqn:iotagh}
\iota(g_n) = t-\iota(u_n) = t+u_n = h_n,
\quad \iota(h_n) = t+\iota(u_n) = t-u_n = g_n.
\end{equation}
Multiplying these two elements, we obtain
\begin{equation}\label{eqn:gnhn}
g_n \cdot h_n = (t-u_n)\cdot (t+u_n) = t^2-u_n^2 = h_n\cdot g_n,
\end{equation}
and since $u_n^2$ is in the center of $D$, 
$g_n\cdot h_n = h_n\cdot g_n$ is an element in the center of $S$.
These calculations are used in the next lemma:

\begin{lem}\label{lem:iotaIn}
For all $n\in \mathbb{Z}$,
the map $\iota:S \to S$ defined in \zcref{eqn:iotaExtended},
satisfies:
\[
\iota(\mathcal{I}_n) = \mathcal{I}_n.
\]
\end{lem}
\begin{proof}
Take $f \in \mathcal{I}_n$.
According to \zcref{prop:idealizer}(1),
this implies that $f \cdot g_n \cdot S \subseteq g_n \cdot S$.
Hence, there exists $s \in S$,
such that $f\cdot g_n = g_n \cdot s$.
Applying $\iota$ to this equality, 
and using the fact $\iota$ reverses multiplication and \zcref{eqn:iotagh}, this gives
\[
h_n \cdot \iota(f) = \iota(s) \cdot h_n.
\]
From which we obtain that:
\[
g_n \cdot (h_n \cdot \iota(f)) \cdot g_n = g_n \cdot (\iota(s) \cdot h_n) \cdot g_n.
\]
Since by \zcref{eqn:gnhn}, $g_n\cdot h_n = h_n\cdot g_n$ is a non-zero central element, and since $S$ is a domain,
we deduce that
\[
\iota(f) \cdot g_n = g_n \cdot \iota(s).
\]
This implies that $\iota(f) \cdot g_n \cdot S \subseteq g_n \cdot S$,
which by \zcref{prop:idealizer}(1) says that $\iota(f) \in \mathcal{I}_n$.
Since $\iota\circ \iota = 1_S$, this implies that $\iota(\mathcal{I}_n) = \mathcal{I}_n$.
\end{proof}

We may now prove the promised symmetry relation for $R$:
\begin{prop}\label{prop:symmetry}
The map $\iota$ is a $\k$-algebra isomorphism $\iota:R \to R^{\opn{op}}$.    
\end{prop}
\begin{proof}
Since $\iota$ is a $\k$-algebra isomorphism $\iota:S \to S^{\opn{op}}$,
it is enough to show that $\iota$ restricts to a bijection $R\to R$.
However, since $R = \cap_{n=1}^{\infty} \mathcal{I}_n$,
we see that
\[
\iota(R) = \iota(\bigcap_{n=1}^{\infty} \mathcal{I}_n) =
\bigcap_{n=1}^{\infty} \iota(\mathcal{I}_n) =
\bigcap_{n=1}^{\infty} \mathcal{I}_n = R,
\]
where the second equality is because $\iota$ is injective, 
and the third equality is \zcref{lem:iotaIn}.
\end{proof}

\subsection{Noetherianity}

Next, we shall prove that $R$ is two-sided noetherian.

\begin{prop}\label{prop:sNoetherian}
The ring $S = D[t]$ is two-sided noetherian.    
\end{prop}
\begin{proof}
Since $D$ is a division ring,
it is two-sided noetherian, so
this is the noncommutative Hilbert basis theorem,
which is a specific case of \cite[Theorem 1.2.9(iv)]{NCBook}.
\end{proof}

Recall that by \zcref{prop:idealizer}, $\mathcal{I}_n$ was shown to be the idealizer of the right ideal $g_n \cdot S$. It turns out, these right ideals are maximal:

\begin{lem}\label{lem:gnMaximal}
For any $n \in \mathbb{Z}$,
the element $g_n = t- u_n \in S$ satisfies that $g_n \cdot S$ is a maximal right ideal of $S$.
\end{lem}
\begin{proof}
It follows from the two parts of \zcref{lem:belongTo} that there is a right $D$-linear decomposition
\[
S = g_n \cdot S \oplus D.
\]
Thus, $S/g_n\cdot S$ is a one dimensional right $D$-module,
and hence, simple.
\end{proof}

The idealizers $\mathcal{I}_n$ are themselves right-noetherian:
\begin{prop}\label{prop:idealizerNoeth}
For any $n \in \mathbb{Z}$,
the idealizer $\mathcal{I}_n$ is a right-noetherian ring.
\end{prop}
\begin{proof}
By \zcref{prop:idealizer},
$\mathcal{I}_n$ is the idealizer of the right ideal $g_n\cdot S$,
and by \zcref{lem:gnMaximal}, this is a maximal right ideal.
Since by \zcref{prop:sNoetherian},
the ring $S$ is right-noetherian,
\cite[Theorem 1.1.12(iii)]{NCBook} says that $\mathcal{I}_n$ is right-noetherian.
\end{proof}

We wish to analyze the right ideals of $R$.
The next lemma is a tool that will soon be used for that goal.

\begin{lem}\label{lem:mult}
Let $F \subseteq \mathbb{Z}$ be a finite subset, 
and let $s \in S$ be some polynomial with coefficients in $D_F$.
For any $m \in \mathbb{Z}$,
consider the element
\[
q_m := \left(\prod_{n \in F, n\ge 1, n \ne m} g_n \cdot h_n\right) \in S.
\]
Then for all $n \ge 1$ such that $n \ne m$,
it holds that $q_m \cdot s\in \mathcal{I}_n$.
If moreover $s \in \mathcal{I}_m$, then $q_m \cdot s \in R$.
In particular, if $m \le 0$ then $q_m \cdot s \in R$.
\end{lem}
\begin{proof}
Let $n \ge 1$,
and suppose that $n \ne m$.
If $n \in F$, using the fact established in \zcref{eqn:gnhn} that $g_n \cdot h_n = h_n \cdot g_n$ is central in $S$,
we may write
\[
q_m \cdot s = g_n \cdot h_n \cdot q_m' \cdot s \in g_n \cdot S \subseteq \mathcal{I}_n,
\]
where $q_m'$ is the element obtained by replacing $F$ with $F\setminus\{n\}$.
If, on the other hand, $n \notin F$,
then as all coefficients of $q_m\cdot s$ lie in $D_F$,
so by \zcref{lem:coeffBelong}, we have that $q_m \cdot s \in \mathcal{I}_n$.
If $s \in \mathcal{I}_m$,
since $q_m$ is central in $S$,
by \zcref{prop:idealizer}(2), $q_m \in \mathcal{I}_m$,
so $q_m\cdot s \in \mathcal{I}_m$.
Thus, in this case, $q_m\cdot s \in \mathcal{I}_n$ for all $n \ge 1$,
so $q_m \cdot s \in R$.
Finally, if $m \le 0$, then any $n \ge 1$ satisfies $n \ne m$,
so we have seen that $q_m \cdot s \in \mathcal{I}_n$ for all $n \ge 1$,
which shows that $q_m \cdot s \in R$.
\end{proof}

We let $Z(D)$ denote the center of $D$,
which is a subfield of the division ring $D$.
Then $Z(D)[t] \subseteq S$ is the center of $S$.
It turns out that every right ideal of $R$ intersects it:

\begin{lem}\label{lem:idealintersectcenter}
Let $J$ be a non-zero right ideal of $R$.
Then $J \cap Z(D)[t] \ne 0$.
In other words, any non-zero right ideal of $R$ contains a non-zero polynomial with coefficients in $Z(D)$.
\end{lem}
\begin{proof}
Take any non-zero element $0 \ne f \in J$.
We shall use several properties of $D$ recorded in \zcref{thm:divisionRing} freely in what follows.
Since $f$ is a finite sum,
we can find a finite subset $F \subseteq \mathbb{Z}$,
such that $f \in D_F[t]$.
Since $D_F$ is a finite dimensional vector space over $Z_F$,
it follows that $D_F[t]$ is a finitely generated free $Z_F[t]$-module.
Let $N$ be its rank.
Since $Z_F[t]$ is contained in the center of $S$,
the map $\lambda:D_F[t] \to D_F[t]$ given by $\lambda(-) = f\cdot -$ is $Z_F[t]$-linear, and since $S$ is a domain, $\lambda$ is injective.
It follows that $c = \det(\lambda) \ne 0$.
We now wish to use the Cayley–Hamilton theorem,
which is valid over the commutative integral domain $Z_F[t]$ on $\lambda$.
If $x^N + \alpha_{N-1}x^{N-1} + \dots + (-1)^N \cdot c$ is the characteristic polynomial of $\lambda$,
we deduce that the element $y = f^{N-1}+\alpha_{N-1}\cdot f^{N-2} + \dots + \alpha_1 \in D_F[t]$ satisfies that $f\cdot y = (-1)^{N+1} \cdot c$.
Let 
\[
q_0 = \prod_{n\in F, n \ge 1} g_n \cdot h_n,
\]
as in \zcref{lem:mult}.
Then by \zcref{lem:mult} with $m=0$,
we have that $q_0\cdot y \in R$.
Hence, since $J$ is a right ideal of $R$, and $f \in J$,
we conclude that $f\cdot q_0\cdot y\in J$.
But $q_0$ is in the center of $R$,
so
\[
f\cdot q_0 \cdot y = q_0 \cdot f\cdot y = q_0 \cdot (-1)^{N+1} \cdot c,
\]
and since $c \in Z_F[t]$ and $q_0 \in Z(D)[t]$ are both non-zero elements and $Z(D)[t]$ is an integral domain,
we obtain the required non-zero element in $J$ which lies in $Z(D)[t]$.
\end{proof}

To prove that $R$ is right-noetherian,
we will do a local-to-global argument over its central subring $Z(D)[t]$. This is a polynomial ring over a field, so has Krull dimension $1$.
In particular, its prime spectrum consists of zero, and the maximal ideals, and these, in turn, are exactly $(f)$ where the $f$'s are all the irreducible monic polynomials over $Z(D)$. One important class of such irreducible polynomials is analyzed in the next result:
 
\begin{lem}\label{lem:irrgnhn}
For every $n \ge 1$,
the element $g_n \cdot h_n = t^2 - u_n^2 \in Z(D)[t]$ is an irreducible monic polynomial. Moreover, if $n \ne m$, $n,m \ge 1$,
then the prime ideals $(g_n \cdot h_n)$ and $(g_m \cdot h_m)$ in $Z(D)[t]$ are not equal.
\end{lem}
\begin{proof}
If $t^2-u_n^2$ was reducible,
then it would have some root $z \in Z(D)$ which would satisfy $z^2 = u_n^2$. Since such a $z$ is a central element, this would imply that $(u_n-z)\cdot(u_n + z) = 0$, and since this is an equation in a domain, 
we would have that $z=\pm u_n$, which is impossible, because $u_n$ is not in the center of $D$ as it does not commute with $v_n$.
Next, let $n,m \ge 1$,
and suppose that $(t^2-u_n^2) = (t^2-u_m^2)$ as prime ideals.
Since both of these are monic polynomials, 
this implies that $u_n^2 = u_m^2$.
Since $u_n$ and $u_m$ commute, 
we deduce that $(u_n-u_m)\cdot(u_n+u_m) = 0$, 
so that $u_n = \pm u_m$.
This is impossible for $n\ne m$,
since then $u_n$ commutes with $v_m$ but not with $v_n$,
while $\pm u_m$ commutes with $v_n$ but not with $v_m$.
Thus, $n=m$.
\end{proof}

We may now calculate the localization of $R$ at every non-zero prime ideal of its central subring $Z(D)[t]$.
By definition, if $T$ is any $Z(D)[t]$-algebra such that $Z(D)[t]$ is contained in its center, and if $\p \in \opn{Spec}(Z(D)[t])$,
we let $T_{\p} := T\otimes_{Z(D)[t]} Z(D)[t]_{\p}$.
If $T\subset U$ are two such $Z(D)[t]$-algebras,
we will identify, using flatness of localization, $T_{\p}$ as a subring of $U_{\p}$.

\begin{lem}\label{lem:localizationofR}
Let $0 \ne \p \in \opn{Spec}(Z(D)[t])$.
Then the following holds:
\begin{enumerate}
\item If $\p \ne (t^2-u_n^2)$ for all $n \ge 1$,
then $R_{\p} = S_{\p}$.
\item If $\p = (t^2-u_n^2)$ for some $n \ge 1$,
then $R_{\p} = \left(\mathcal{I}_n\right)_{\p}$.
\end{enumerate}
\end{lem}
\begin{proof}
\quad
\begin{enumerate}
\item Since $R \subseteq S$, 
we know that $R_{\p} \subseteq S_{\p}$.
To see the converse, it is enough to show that $S \subseteq R_{\p}$.
Take any $0 \ne s \in S$,
and take some finite subset $F \subseteq \mathbb{Z}$,
such that $s \in D_F[t]$.
As in \zcref{lem:mult},
let 
\[
q_0 = \prod_{n \in F, n\ge 1} g_n\cdot h_n.
\]
Applying \zcref{lem:mult} with $m=0$, 
we see that $q_0\cdot s \in R$.
We know that for any $1 \le n \in F$,
$g_n\cdot h_n = t^2-u_n^2 \notin \p$,
because \zcref{lem:irrgnhn} showed that $(t^2-u_n^2)$ is a maximal ideal,
and $\p$ is a different maximal ideal.
Hence, the fact that $\p$ is prime implies that $q_0 \notin \p$.
We see that $s = q_0^{-1}\cdot q_0 \cdot s \in R_{\p}$.
\item Assume that $\p = (t^2-u_n^2)$ for some $n \ge 1$.
Since $R \subseteq \mathcal{I}_n$,
we know that $R_{\p} \subseteq \left(\mathcal{I}_n\right)_{\p}$.
As before, it is enough to show that $\mathcal{I}_n \subseteq R_{\p}$.
Let $s \in \mathcal{I}_n$,
and let $F \subseteq \mathbb{Z}$ be a finite set such that $s \in D_F[t]$.
Take
\[
q_n = \prod_{k \in F, k\ge 1, k\ne n} g_k\cdot h_k.
\]
Since $s \in \mathcal{I}_n$, by \zcref{lem:mult},
we know that $q_n\cdot s \in R$.
According to \zcref{lem:irrgnhn},
we know that for $k \ne n$, $g_k \cdot h_k \notin \p$,
because $\p$ and $(g_k\cdot h_k)$ are two different maximal ideals.
Hence, the fact that $\p$ is prime implies that $q_n \notin \p$,
which shows that $s = (q_n)^{-1} \cdot q_n \cdot s \in R_{\p}$.
\end{enumerate}
\end{proof}

We may now prove that $R$ is locally right-noetherian:

\begin{prop}\label{prop:locals-are-noetherian}
For every maximal ideal $0 \ne \p \in \opn{Spec}(Z(D)[t])$,
the ring $R_{\p}$ is right-noetherian.
\end{prop}
\begin{proof}
We consider the two cases as in \zcref{lem:localizationofR}.
Suppose first that $\p \ne (t^2-u_n^2)$ for all $n \ge 1$,
so by \zcref{lem:localizationofR}, we know that $R_{\p} = S_{\p}$.
By \zcref{prop:sNoetherian} we know that $S$ is right-noetherian.
According to \cite[Proposition 2.1.16(iii)]{NCBook},
every right ideal $J$ of $R_{\p} = S_{\p}$ is of the form $J = (J\cap S)\cdot S_{\p}$, which is clearly finitely generated since $J \cap S$ is.
Next, assume that $\p = (t^2-u_n^2)$ for some $n \ge 1$.
Then \zcref{lem:localizationofR} tells us that $R_{\p} = \left(\mathcal{I}_n\right)_{\p}$.
Since by \zcref{prop:idealizerNoeth}, the ring $\mathcal{I}_n$ is right-noetherian, repeating the same argument from the first part of this proof shows that $R_{\p}$ is right-noetherian in this case as well.
\end{proof}

Next, we prove a general local-to-global result for noetherianity that suits our situation:
\begin{prop}\label{prop:local-to-global-noetherian}
Let $B$ be a commutative principal ideal domain which is not a field,
let $C$ be a ring, and assume that $B$ is contained in the center of $C$.
Suppose that for any maximal ideal $\m$ of $B$,
it holds that $C_{\m} := C\otimes_B B_{\m}$ is a right-noetherian ring.
Suppose further that for any non-zero right ideal $J$ of $C$,
we have that $J \cap B \ne 0$.
Then $C$ is right-noetherian.
\end{prop}
\begin{proof}
Let $J$ be a non-zero right ideal of $C$.
We must show that $J$ is finitely generated.
By assumption, there is a non-zero element $b \in J$,
such that $b \in B$.
If $b$ is a unit, then $J = (1)$, and we are done.
Otherwise, using the fact that $B$ is a unique factorization domain,
factor $b$ as a product of different primes
\[
b = \prod_{i=1}^n p_i^{e_i}, \quad e_i \ge 1, \quad i\ne j\implies (p_i) \ne (p_j).
\]
Let us write for any $1 \le i \le n$, $\m_i = (p_i)$.
Note that the $\m_i$ are distinct maximal ideals,
so by the Chinese remainder theorem (\cite[Proposition 1.10]{AM}) over $B$, there is an isomorphism of $B$-modules:
\[
B/b \cong \prod_{i=1}^n B/(p_i^{e_i}).
\]
Applying $-\otimes_B C$,
we obtain an isomorphism of right $C$-modules
\begin{equation}\label{eqn:Cmodb}
C/(b\cdot C) \cong \prod_{i=1}^n C/(p_i^{e_i} \cdot C).
\end{equation}
For every $1 \le i \le n$,
note that $B/(p_i^{e_i})$ is a $B_{\m_i}$-module.
Indeed, if $x \notin \m_i$, then $\gcd(x,p_i^{e_i}) = 1$,
so there exist $u, v \in B$,
such that $u\cdot x + v\cdot p_i^{e_i} = 1$,
which shows that $x$ is invertible in $B/(p_i^{e_i})$.
This implies that 
\[
C/(p_i^{e_i} \cdot C) = C_{\m_i}/(p_i^{e_i} \cdot C_{\m_i})
\]
is a right module over $C_{\m_i} = B_{\m_i}\otimes_B C$,
and clearly it is finitely generated over it.
Since by assumption $C_{\m_i}$ is right-noetherian,
we deduce that $C/(p_i^{e_i} \cdot C)$ is a noetherian right $C_{\m_i}$-module.
But sub $C$-modules and sub $C_{\m_i}$-modules of $C/(p_i^{e_i}\cdot C)$ coincide, because the element $u$ above acts on them as $x^{-1}$,
so we deduce that $C/(p_i^{e_i} \cdot C)$ is a noetherian right $C$-module.
Hence, \zcref{eqn:Cmodb} shows that $C/(b\cdot C)$ is also noetherian as a right $C$-module.
This implies that its submodule $J/(b\cdot C)$ is a finitely generated right $C$-module.
The exact sequence of right $C$-modules
\[
0 \to b\cdot C \to J \to J/(b\cdot C) \to 0
\]
then implies that $J$ is finitely generated as a right module over $C$.
Hence, $C$ is right-noetherian.
\end{proof}

We may finally show that $R$ is noetherian:

\begin{thm}\label{thm:noetherian}
The ring $R$ constructed in \zcref{eqn:theR} is a two-sided noetherian ring.
\end{thm}
\begin{proof}
First, we show that $R$ is right-noetherian.
To do that, let us check that the ring $R$ satisfies the conditions of \zcref{prop:local-to-global-noetherian}.
The principal ideal domain $B$, which is not a field, in that proposition will be $B = Z(D)[t]$,
which we know is a subring of $R$ contained in its center.
The fact that all localizations of $R$ at maximal ideals of $B$ are right-noetherian is \zcref{prop:locals-are-noetherian},
and the fact that any non-zero right ideal of $R$ intersects with $B$ non-trivially is \zcref{lem:idealintersectcenter}.
Hence, we deduce that $R$ is right-noetherian.
According to \zcref{prop:symmetry},
there is a ring isomorphism $R \cong R^{\opn{op}}$,
from which we deduce that $R$ is also left-noetherian.
\end{proof}

\subsection{The hereditary property}

The final goal in this section is to show that $R$ is a hereditary ring,
that is, that it has global dimension at most $1$.
Since we already proved that $R$ is two-sided noetherian,
by Auslander's theorem \cite[Corollary 5]{Aus}, the global dimension of $R$ is left-right symmetric, so it is sufficient to show that $R$ is right hereditary.
The basic strategy will be as in the noetherianity section:
we will show that the result holds locally over the maximal spectrum of $Z(D)[t]$, and then prove a local-to-global result.
This strategy has the advantage that we have already calculated, in \zcref{lem:localizationofR}, all the localizations of $R$ at these maximal ideals.

\begin{prop}\label{prop:shere}
The ring $S$ is right hereditary.
\end{prop}
\begin{proof}
The ring $S$ is a polynomial ring over the division ring $D$.
$D$, being a division ring, has right global dimension $0$,
so by \cite[Theorem 7.5.3(iii)]{NCBook} (applied with $\sigma = 1_D$),
we see that $S$ has right global dimension $1$.
\end{proof}

Next, we show that the idealizers are right hereditary.

\begin{prop}\label{prop:ihere}
For any $n \in \mathbb{Z}$,
the idealizer $\mathcal{I}_n$ is a right hereditary ring.
\end{prop}
\begin{proof}
Let $J = (t-u_n)\cdot S$.
By \zcref{lem:gnMaximal}, $J$ is a maximal right ideal of $S$.
We have just seen in \zcref{prop:shere} that $S$ is right hereditary.
According to \cite[Theorem 5.5.8]{NCBook},
this implies that the idealizer of $S$ with respect to any isomaximal right ideal is also right hereditary.
The notion of an isomaximal right ideal is defined in \cite[5.5.2]{NCBook}, where it is explained that maximal right ideals are isomaximal.
Thus, we deduce that the idealizer of $S$ with respect to $J$ is right hereditary, and by \zcref{prop:idealizer}, this is exactly $\mathcal{I}_n$.
\end{proof}

We may now show that localizations of $R$ at maximal ideals of $Z(D)[t]$ are right hereditary:
\begin{prop}\label{prop:locals-are-here}
For every maximal ideal $0 \ne \p \in \opn{Spec}(Z(D)[t])$,
the ring $R_{\p}$ is right hereditary.    
\end{prop}
\begin{proof}
According to \zcref{lem:localizationofR},
the ring $R_{\p}$ is the central localization of either $S$ or $\mathcal{I}_n$. 
Both $S$ and $\mathcal{I}_n$ are right hereditary by \zcref{prop:shere} and \zcref{prop:ihere}.
Hence, the result follows from \cite[Corollary 7.4.3]{NCBook},
which tells us that central localization cannot increase the right global dimension.
\end{proof}

Now, the local-to-global result:

\begin{prop}\label{prop:local-to-global-here}
Let $B$ be a commutative ring,
let $C$ be a right noetherian ring, and assume that $B$ is contained in the center of $C$.
Suppose that for any maximal ideal $\m$ of $B$,
it holds that $C_{\m} := C\otimes_B B_{\m}$ is right hereditary.
Then $C$ is right hereditary.
\end{prop}
\begin{proof}
Let $J$ be a right ideal of $C$.
According to \cite[Proposition 5.14]{Lam2},
it is enough to show that $J$ is projective over $C$. Equivalently, that $\opn{Ext}^1_C(J,-) = 0$.
Let $M$ be an arbitrary right $C$-module.
The argument now is a standard localization argument,
but since it is difficult to find a reference in such a noncommutative context, we give the full argument.
To compute this $\opn{Ext}^1(J,M)$,
we take a projective resolution of $J$ over $C$.
The fact that $C$ is right-noetherian and $J$ is finitely generated implies we may find such a resolution
\[
\dots \to P^{-1} \to P^0 \to J \to 0
\]
such that each $P^i$ is a finitely generated projective right $C$-module.
Applying the functor $\opn{Hom}_C(-,M)$ to the resolution gives a complex
\[
0 \to \opn{Hom}_C(P^0,M) \to \opn{Hom}_C(P^{-1},M) \to \dots
\]
of $B$-modules.
Take a maximal ideal $\m \subseteq B$.
For each $i \le 0$, the fact that $P^i$ is finitely generated implies that there is a natural isomorphism
\[
\opn{Hom}_C(P^i,M) \otimes_B B_{\m} \cong \opn{Hom}_{C_{\m}}(P^i_{\m},M_{\m}).
\]
Indeed, if $P^i = C^n$ for some $n$ this is clear,
and the general case follows from the fact that both of these functors are additive and $P^i$ is a direct summand of some $C^n$.
Passing to cohomology, and using the fact that $B_{\m}$ is flat over $B$,
we thus deduce that the $B$-module $\opn{Ext}^1_C(J,M)$ satisfies:
\[
\left(\opn{Ext}^1_C(J,M)\right)_{\m} = \opn{Ext}^1_{C_{\m}}(J_{\m},M_{\m}).
\]
By assumption, the ring $C_{\m}$ is right hereditary,
so its right ideal $J_{\m}$ is projective over it.
Hence, $\opn{Ext}^1_{C_{\m}}(J_{\m},M_{\m}) = 0$.
Thus, all localizations of the $B$-module $\opn{Ext}^1_C(J,M)$ at maximal ideals vanish, so by \cite[Proposition 3.8]{AM}, 
we deduce that $\opn{Ext}^1_C(J,M) = 0$.
Hence, $J$ is projective over $C$, which implies that $C$ is right hereditary.
\end{proof}

Finally, $R$ itself is a hereditary ring:
\begin{thm}\label{thm:here}
The ring $R$ constructed in \zcref{eqn:theR} is a hereditary ring. 
\end{thm}
\begin{proof}
The fact that $R$ is right hereditary follows from \zcref{prop:local-to-global-here} and \zcref{prop:locals-are-here},
by taking $B = Z(D)[t]$.
By \zcref{thm:noetherian}, $R$ is noetherian,
so by Auslander's theorem, \cite[Corollary 5]{Aus}, its left and right global dimension coincide,
so $R$ is also left hereditary.
\end{proof}

\subsection{The main result}

All the pieces are in place, and we may finally prove the main result of this section:

\begin{thm}\label{thm:mainR}
Let $\k$ be any field with $\opn{char}(\k) \ne 2$.    
Then there exists a $\k$-algebra $R$ which satisfies the following conditions:
\begin{enumerate}
\item $R$ is a domain which is not a division ring.
\item $R$ is a two-sided noetherian ring.
\item $R$ is hereditary.
\item There is an injective $\k$-algebra homomorphism $\sigma:R \to R$,
which is not surjective, such that $1 \in (R:\sigma(R))_r\cdot R$.
\end{enumerate}
\end{thm}
\begin{proof}
Take $R$ to be the ring constructed in \zcref{eqn:theR}.
By \zcref{prop:rbasic}, $R$ is a $\k$-algebra which is a domain and not a division ring.
By \zcref{thm:noetherian}, $R$ is a two-sided noetherian ring,
and by \zcref{thm:here}, $R$ is left and right hereditary.
Take $\sigma:R \to R$ to be the $\k$-algebra homomorphism defined in \zcref{eqn:sigma}. 
By \zcref{prop:injnotsurj}, $\sigma$ is injective but not surjective,
and the final condition $1 \in (R:\sigma(R))_r\cdot R$ is shown in \zcref{prop:conductor}.
\end{proof}

\begin{cor}\label{cor:main}
Let $\k$ be any field with $\opn{char}(\k) \ne 2$.    
Then there exists a $\k$-algebra $A$ which satisfies:
\begin{enumerate}
\item $A$ is a two-sided noetherian ring.
\item There are equalities
\[
\injdim_A({}_A A) = 1, \quad \injdim_A(A_A) = +\infty.
\]
\end{enumerate}
\end{cor}
\begin{proof}
Take a $\k$-algebra $R$ and a $\k$-algebra map $\sigma:R \to R$ as in \zcref{thm:mainR}.
Let $A = R[x;\sigma]/(x^2)$ as in \zcref{eqn:hero}.
The fact that $\sigma$ is $\k$-linear implies that $A$ is a $\k$-algebra,
and the result now follows from \zcref{thm:main-gen-const} and \zcref{thm:mainR}.
\end{proof}

\begin{cor}
The algebra $A$ in \zcref{cor:main} is not two-sided artinian.    
\end{cor}
\begin{proof}
By \cite[Corollary]{Ho}, a two-sided artinian ring $A$ with $\injdim_A({}_A A) \le 1$ has $\injdim_A(A_A) \le 1$. 
Since $\injdim_A(A_A) = +\infty$, $A$ cannot be two-sided artinian.
\end{proof}

\begin{rem}
One can show directly that $A$ is neither left nor right artinian,
but this won't be needed, so this proof is omitted.
\end{rem}

\section{Dualizing complexes}\label{sec:dual}

Dualizing complexes were first introduced by Grothendieck and recorded by Hartshorne in \cite{RD} in commutative algebra and algebraic geometry contexts.
They were later introduced to noncommutative algebra by Yekutieli \cite{Yek},
and substantively developed in that context by Yekutieli-Zhang \cite{YZ1,YZ2}, Van den Bergh \cite{vdb}, Iyengar-Krause \cite{InKr}, Neeman \cite{Neeflat} and Positselski \cite{Pos}. The goal of this section is to show that the algebra $A$ from \zcref{cor:main} has a dualizing complex.
This will be used in the next section to obtain further results about $A$.

Dualizing complexes are objects of derived categories. We refer the reader to \cite{KelBook,YeBook} for basics about derived categories.
For a ring $B$, we denote by $\cat{D}({}_B \opn{Mod})$ the unbounded derived category of left $B$-modules.
In particular, if $B$ is a $\k$-algebra over a field $\k$,
then objects of $\cat{D}({}_{B\otimes_{\k} B^{\opn{op}}} \opn{Mod})$ are complexes of $B$-$B$-bimodules.

We define the restriction functors 
\begin{equation}\label{eqn:Restf}
\restlb: {}_{B\otimes_{\k} B^{\opn{op}}}\opn{Mod} \to {}_B \opn{Mod}, \quad
\restrb: {}_{B\otimes_{\k} B^{\opn{op}}}\opn{Mod} \to \opn{Mod}_B
\end{equation}
to be the functors that forget the right and left structures of a bimodule respectively.
These naturally extend to complexes of bimodules, as well as to $\cat{D}({}_{B\otimes_{\k} B^{\opn{op}}} \opn{Mod})$ (since they are exact).

For a $\k$-algebra $B$, assumed to be two-sided noetherian, we define a dualizing complex $\omega$ to be an object of $\cat{D}({}_{B\otimes_{\k} B^{\opn{op}}} \opn{Mod})$ which satisfies the following conditions:
\begin{enumerate}
\item The restrictions $\restlb(\omega)$ and $\restrb(\omega)$ have bounded cohomology.
\item For all $n \in \mathbb{Z}$, the left $B$-module $\mrm{H}^n(\restlb(\omega))$ is a finitely generated left $B$-module.
\item For all $n \in \mathbb{Z}$, the right $B$-module $\mrm{H}^n(\restrb(\omega))$ is a finitely generated right $B$-module.
\item It holds that $\injdim_B(\restlb(\omega)) < \infty$,
and that $\injdim_B(\restrb(\omega))<\infty$.
\item The homothety maps 
\[
B \to \mrm{R}\opn{Hom}_{B^{\opn{op}}}(\omega,\omega), \quad
B^{\opn{op}} \to \mrm{R}\opn{Hom}_B(\omega,\omega)
\]
are isomorphisms in $\cat{D}({}_{B\otimes_{\k} B^{\opn{op}}} \opn{Mod})$.
\end{enumerate}

The point of this definition is that if $\omega$ is a dualizing complex over the $\k$-algebra $B$,
then by \cite[Proposition 1.3]{YZ1}, for any complex $M$ of left $B$-modules, with $\mrm{H}^n(M)$ being finitely generated for all $n$,
the derived duality map
\[
M \to \mrm{R}\opn{Hom}_{B^{\opn{op}}}(\mrm{R}\opn{Hom}_B(M,\omega),\omega)
\]
is an isomorphism in $\cat{D}({}_B \opn{Mod})$.
A similar result holds for complexes of right $B$-modules.

Before we can prove that the algebra $A$ from \zcref{cor:main} has a dualizing complex, we must gather additional facts about the ring $R$ constructed in \zcref{thm:mainR}.

\subsection{The conductor $(R:\sigma(R))_r$}

This subsection can be seen as a continuation of \zcref{sec:conductor}.
In particular, we will freely use the notation set up in that section.
Our goal is to show that the left ideal 
\[
(R:\sigma(R))_r = \{x \in R \mid R\cdot x \subseteq \sigma(R)\}
\]
is the principal left ideal generated by 
\begin{equation}\label{eqn:pi}
\pi_0 = g_0\cdot v_0 = (t-u_0)\cdot v_0 = t\cdot v_0 -u_0 \cdot v_0.
\end{equation}

The maps $\beta, \iota:S \to S$ were defined in \zcref{sec:sigmaDef} and \zcref{sec:symm} respectively.

\begin{lem}\label{lem:iotabeta}
There is an equality $\iota\circ \beta = \beta \circ \iota$ as maps $S \to S$.
\end{lem}
\begin{proof}
We may view $\iota\circ \beta$ and $\beta\circ \iota$ as $\k$-algebra maps $S \to S^{\opn{op}}$.
Let 
\[
T = \{x \in S\mid \iota\circ \beta(x) = \beta \circ \iota(x)\}.
\]
Clearly $T$ is a sub $\k$-algebra of $S$.
For any $n \in \mathbb{Z}$,
we have
\[
\iota(\beta(u_n)) = \iota(u_{n-1}) = -u_{n-1},
\quad
\beta(\iota(u_n)) = \beta(-u_n) = -u_{n-1},
\]
and similarly 
\[
\iota(\beta(v_n)) = \iota(v_{n-1}) = v_{n-1},
\quad
\beta(\iota(v_n)) = \beta(v_n) = v_{n-1}.
\]
Thus, for all $n \in \mathbb{Z}$, $u_n, v_n \in T$.
Since both $\iota \circ \beta$ and $\beta \circ \iota$ reverse products,
they respect inverses, so $T \cap D$ is a division ring.
Since $D$ is generated by the $u_n, v_n$ as a division ring,
we deduce that $D \subseteq T$.
The construction of $\beta, \iota$ further says that:
\[
\iota(\beta(t)) = t = \beta(\iota(t)),
\]
so $t \in T$. Since $S$ is generated as an algebra by $D$ and $t$,
we deduce that $T = S$.
\end{proof}

\begin{prop}\label{prop:ellmult}
Let $f,g \in S$,
and let $d \in D$.
Suppose that $f_{\ell}(d)\cdot d = d\cdot f_{\ell}(d)$.
Then 
\[
(f\cdot g)_{\ell}(d) = f_{\ell}(d) \cdot g_{\ell}(d).
\]
\end{prop}
\begin{proof}
Let $c = f_{\ell}(d)$.
Suppose first $c \ne 0$.
Then by the left version of \cite[Proposition 16.3]{Lam},
we have
\[
(f\cdot g)_{\ell}(d) = f_{\ell}(d) \cdot g_{\ell}(c^{-1}\cdot d \cdot c).
\]
Since by assumption $c\cdot d = d\cdot c$, we deduce that $c^{-1}\cdot d \cdot c = d$, which gives the result.
On the other hand, if $c = 0$,
by the left version of \cite[Proposition 16.2]{Lam},
$f \in (t-d)\cdot S$, so $f\cdot g \in (t-d)\cdot S$,
so using \cite[Proposition 16.2]{Lam} again, we see that $(f\cdot g)_{\ell}(d)$ is also $0$.
\end{proof}

Extending \zcref{eqn:pi},
for $n \in \mathbb{Z}$, let 
\begin{equation}\label{eqn:pin-def}
\pi_n := g_n \cdot v_n = (t-u_n)\cdot v_n \in S.
\end{equation}

Here are some elementary computations concerning this element:
\begin{lem}\label{lem:piIdeal}
The element $\pi_n$ defined above satisfies:
\begin{enumerate}
\item There is an equality: $\pi_n = v_n \cdot h_n = v_n \cdot (t+u_n)$.
\item The right ideals $g_n\cdot S = (t-u_n)\cdot S$ and $\pi_n \cdot S$ coincide.
\item For all $n \ge 0$, the element $\pi_n$ belongs to $R$.
\end{enumerate}
\end{lem}
\begin{proof}\quad
\begin{enumerate}
\item Indeed, since $t$ is in the center of $S$,
and $u_n,v_n$ anti-commute:
\[
(t-u_n)\cdot v_n = t\cdot v_n -u_n\cdot v_n = v_n\cdot t + v_n\cdot u_n = v_n\cdot(t+u_n).
\]
\item The right ideals $g_n \cdot S$ and $g_n \cdot v_n \cdot S$ coincide because $v_n$ is a unit in $S$.
\item If $n \ge 1$, then $\pi_n \in D_{\{n\}}[t]$,
so by \zcref{lem:coeffBelong}, $\pi_n \in \mathcal{I}_m$ for all $m \ne n$. By \zcref{prop:idealizer}(2), $g_n\cdot S \subseteq \mathcal{I}_n$,
so (2) shows that $\pi_n \in \mathcal{I}_n$.
Thus, by definition of $R$, $\pi_n \in R$.
Finally, if $n = 0$, then $\pi_0 \in D_{\{0\}}[t]$ and \zcref{lem:coeffBelong} implies $\pi_0 \in \mathcal{I}_n$ for all $n\ge 1$, so $\pi_0 \in R$.
\end{enumerate}
\end{proof}

The special elements $\pi_m$ satisfy a cancellation property:

\begin{lem}\label{lem:canc}
For any $m \in \mathbb{Z}$,
given $s \in S$,
if $\pi_m\cdot s \in R$, then $s \in \mathcal{I}_n$ for all $n \ge 1$ such that $n \ne m$.
If moreover $m \le 0$, then $s \in R$.
\end{lem}
\begin{proof}
By definition of $\pi_m$, for all $n \ge 1$ with $n \ne m$, we have:
\[
(\pi_m)_{\ell}(u_n) = (t\cdot v_m-u_m \cdot v_m)_{\ell}(u_n) = 
u_n\cdot v_m-u_m\cdot v_m = (u_n-u_m)\cdot v_m.
\]
This element is non-zero because $u_n$ commutes with $v_m$ but $u_m$ does not.
Moreover, it commutes with $u_n$ because $u_m$ and $v_m$ commute with $u_n$, as $n \ne m$.
Let $f = \pi_m$, $g=s$, and $d = u_n$.
Then by \zcref{prop:ellmult}:
\begin{equation}\label{eqn:multell}
(\pi_m\cdot s)_{\ell}(u_n) = (\pi_m)_{\ell}(u_n) \cdot s_{\ell}(u_n).
\end{equation}
The assumption that $\pi_m\cdot s \in R$ implies in particular that $\pi_m\cdot s \in \mathcal{I}_n$,
so the definition of $\mathcal{I}_n$ tells us that
\[
(\pi_m\cdot s)_{\ell}(u_n) \in C_D(u_n).
\]
We have seen that $(\pi_m)_{\ell}(u_n)$ is a non-zero element of the division ring $C_D(u_n)$,
so we deduce from \zcref{eqn:multell} that
\[
s_{\ell}(u_n) = \left((\pi_m)_{\ell}(u_n)\right)^{-1} \cdot (\pi_m\cdot s)_{\ell}(u_n)
\]
also belongs to $C_D(u_n)$.
By definition of $\mathcal{I}_n$,
this implies that $s \in \mathcal{I}_n$.
Since this holds for all $n \ge 1$ with $n \ne m$, the first claim follows. If $m \le 0$, every $n \ge 1$ satisfies $n \ne m$,
which implies that $s \in \mathcal{I}_n$ for all $n \ge 1$.
Hence, $s \in R$.
\end{proof}

The above is enough to compute the left conductor of $(R:\sigma(R))$:

\begin{prop}\label{prop:left-conductor}
For the element $\pi_0$ defined in \zcref{eqn:pi},
there is an equality of right ideals:
\[
\{r \in R \mid r\cdot R \subseteq \sigma(R)\} = \pi_0 \cdot R.
\]
\end{prop}
\begin{proof}
Note that by \zcref{lem:piIdeal}
(2, 3),
\[
\pi_0 \cdot R \subseteq R \cap \pi_0 \cdot S = R \cap g_0 \cdot S.
\]
By \zcref{prop:idealizer}(2), we know that $g_0 \cdot S \subseteq \mathcal{I}_0$,
so we obtain
\[
\pi_0 \cdot R \subseteq R \cap \mathcal{I}_0 = \sigma(R),
\]
where the equality is \zcref{prop:imsigma}.
Thus, if $r \in R$ satisfies $r \in \pi_0 \cdot R$,
then
\[
r\cdot R \subseteq \pi_0 \cdot R \subseteq \sigma(R),
\]
which gives
\[
\pi_0 \cdot R \subseteq \{r \in R \mid r\cdot R \subseteq \sigma(R)\}.
\]
Conversely, take $r \in R$ such that $r\cdot R \subseteq \sigma(R)$.
In particular, as $1, v_0 \in R$, we deduce that $r, r\cdot v_0 \in \sigma(R) = R \cap \mathcal{I}_0$.
Thus, $r, r\cdot v_0 \in \mathcal{I}_0$.
Let $e = r_{\ell}(u_0)$.
The fact that $r \in \mathcal{I}_0$ implies that $e$ commutes with $u_0$.
We now apply \zcref{prop:ellmult} with $f = r$, $g = v_0$ and $d = u_0$.
This gives us:
\[
(r\cdot v_0)_{\ell}(u_0) = r_{\ell}(u_0) \cdot v_0 = e \cdot v_0.
\]
But $r\cdot v_0$ also belongs to $\mathcal{I}_0$,
so we deduce that $e\cdot v_0$ commutes with $u_0$.
On the other hand, a direct computation gives:
\[
(e\cdot v_0) \cdot u_0 =  - e\cdot u_0 \cdot v_0 = - u_0 \cdot (e\cdot v_0).
\]
Since $\opn{char}(\k) \ne 2$,
this implies that $e\cdot (v_0 \cdot u_0) = 0$,
and since $D$ is a domain and $v_0 \cdot u_0 \ne 0$,
we deduce that $e = 0$.
By \zcref{lem:belongTo}, the fact that $e = r_{\ell}(u_0) = 0$ implies that $r \in g_0 \cdot S = \pi_0 \cdot S$ where the equality of right ideals is \zcref{lem:piIdeal}(2).
We obtained that there is some $s \in S$ such that $r = \pi_0 \cdot s$.
As $\pi_0 \cdot s = r \in R$, by the cancellation property of \zcref{lem:canc}, we see that $s \in R$.
This shows that $r = \pi_0 \cdot s \in \pi_0\cdot R$, 
as claimed.
\end{proof}

We use this to compute the right conductor,
the main result of this subsection:

\begin{prop}\label{prop:rightcondComp}
For the element $\pi_0$ defined in \zcref{eqn:pi},
there is an equality of left ideals:
\[
\{r \in R \mid R\cdot r \subseteq \sigma(R)\} = R \cdot \pi_0.
\]
\end{prop}
\begin{proof}
Consider the equality
\begin{equation}\label{eqn:leftcondeq}
\{r \in R \mid r\cdot R \subseteq \sigma(R)\} = \pi_0 \cdot R,
\end{equation}
proved in \zcref{prop:left-conductor},
and apply the map $\iota$ to both sides of it.
Note first that 
\[
\iota(\pi_0) = \iota(t\cdot v_0 -u_0\cdot v_0) = 
t\cdot v_0 - u_0 \cdot v_0 = \pi_0.
\]
Since $\iota$ is additive and reverses multiplication,
we deduce that 
\begin{equation}\label{eqn:iotapileftright}
\iota(\pi_0 \cdot R) = \iota(R) \cdot \iota(\pi_0) = R\cdot \pi_0,
\end{equation}
where the last equality is because $\iota:R \to R$ is surjective.
Take $r \in R$ such that $R\cdot r \subseteq \sigma(R)$.
Applying $\iota$, we obtain that
\[
\iota(r) \cdot R \subseteq \iota(\sigma(R)) = \sigma(\iota(R)) = \sigma(R).
\]
Here, the first inclusion is because $\iota$ reverses multiplication,
the second equality is \zcref{lem:iotabeta} and the fact that $\sigma$ is the restriction of $\beta$, and the last because $\iota:R \to R$ is surjective.
Conversely, if $\iota(r)\cdot R \subseteq \sigma(R)$,
applying $\iota$, and using the fact that $\iota\circ \iota = 1_R$,
we see that 
\[
R\cdot r \subseteq \iota(\sigma(R)) = \sigma(R).
\]
This shows that $\iota$ maps $\{r \in R \mid r\cdot R \subseteq \sigma(R)\}$ onto $\{r \in R \mid R\cdot r \subseteq \sigma(R)\}$. 
This, \zcref{eqn:leftcondeq} and \zcref{eqn:iotapileftright} give the result.
\end{proof}

\subsection{The automorphism $\tau:R \to R$}

Our next goal is to construct a $\k$-linear automorphism $\tau:R\to R$,
with the special property that $\tau(\sigma(r)) \cdot \pi_0 = \pi_0 \cdot \sigma(\tau(r))$ for all $r \in R$. This will be the final input we will need about $R$, to construct a dualizing complex over the algebra $A = R[x;\sigma]/(x^2)$.

For $n \ge 0$, The elements $\pi_n \in R$ were defined in \zcref{eqn:pin-def}.
For any $N \ge 1$, we now define a new element:
\begin{equation}\label{eqn:PN}
P_N := \pi_1 \cdot \dots \cdot \pi_N \in R,
\end{equation}
and set also $P_0 = 1$.
Here are some elementary properties of $\pi_n$ and $P_N$ we will require:

\begin{lem}\label{lem:elempin}
The following holds for all $n \ge 0$:
\begin{enumerate}
\item It holds that $\beta(\pi_n) = \pi_{n-1}$.
\item If $F \subseteq \mathbb{Z}$ is any finite set such that $n\notin F$, then $\pi_n$ commutes with every element of $D_F[t]$.
\item For $n \ne m$, $n,m \ge 0$, the elements $\pi_n, \pi_m$ commute.
\item If $n\ge 1$, it holds that $\beta(P_n) = \pi_0 \cdot P_{n-1}$, and $\pi_0$ commutes with $P_{n-1}$.
\end{enumerate}
\end{lem}
\begin{proof}\quad
\begin{enumerate}
\item By definition,
\[
\beta(\pi_n) = \beta(t\cdot v_n - u_n\cdot v_n) = 
t\cdot \beta(v_n) - \beta(u_n)\cdot\beta(v_n) = 
t\cdot v_{n-1} - u_{n-1}\cdot v_{n-1} = \pi_{n-1}.
\]
\item By \zcref{thm:divisionRing}(3)(c),
we know that $u_n$ and $v_n$ each commutes with every element of $D_F[t]$, and $t$ is in the center of $D_F[t]$,
so we deduce that $\pi_n = t\cdot v_n - u_n\cdot v_n$ also commutes with every element of $D_F[t]$.
\item Taking $F = \{m\}$, we have that $\pi_m \in D_F[t]$,
so (3) follows from (2).
\item Using the multiplicativity of $\beta$ and $(1)$, calculate:
\[
\beta(P_n) = \beta(\pi_1\cdot \dots \cdot \pi_n) = 
\beta(\pi_1) \cdot \dots \cdot \beta(\pi_n) = 
\pi_0 \cdot \pi_1 \cdot \dots \cdot \pi_{n-1} = \pi_0 \cdot P_{n-1}.
\]
The fact that $\pi_0$ commutes with $P_{n-1}$ follows from (3).
\end{enumerate}
\end{proof}

The elements $g_n, h_n$ satisfy a swap property which will be useful soon:
\begin{lem}\label{lem:swap}
For any $n \in \mathbb{Z}$ and any $s \in S$,
it holds that $s\cdot g_n \in g_n\cdot S$ if and only if $h_n \cdot s \in S \cdot h_n$.
\end{lem}
\begin{proof}
By \zcref{eqn:gnhn} we know that $g_n \cdot h_n = h_n \cdot g_n$ is central in $S$.
Suppose first that $s\cdot g_n \in g_n\cdot S$.
Then $s\cdot g_n = g_n \cdot s'$ for some $s' \in S$.
Multiplying this by $h_n$ on the left, and using the centrality of $g_n\cdot h_n$,
we may cancel $g_n$ and obtain that $h_n \cdot s = s' \cdot h_n$,
so $h_n \cdot s \in S\cdot h_n$.
The proof of the converse is identical with the roles of $g_n$ and $h_n$ reversed.
\end{proof}

Next, we identify the one-sided ideals that $\pi_n$ generate.
This requires some easy but long computations:

\begin{lem}\label{lem:pi-left-right-ideals}
For any $n \ge 1$, 
there are equalities:
\[
\pi_n \cdot R = R \cap g_n \cdot S \cap S\cdot h_n = R\cdot \pi_n.
\]
\end{lem}
\begin{proof}
Let us set $J:= R \cap g_n \cdot S \cap S\cdot h_n$.
Step 1: $\pi_n\cdot R \subseteq J$:
Since $\pi_n \in R$, we know that $\pi_n\cdot R \subseteq R$.
Secondly, $R \subseteq S$ implies $\pi_n\cdot R \subseteq \pi_n \cdot S = g_n \cdot S$,
where the last equality is \zcref{lem:piIdeal}(2).
Finally, we show that $\pi_n \cdot R \subseteq S\cdot h_n$:
take $r \in R$. In particular, $r \in \mathcal{I}_n$.
By \zcref{lem:swap}, $h_n\cdot r \in S\cdot h_n$. Since $\pi_n = v_n\cdot h_n$ by \zcref{lem:piIdeal}(1), we get $\pi_n\cdot r = v_n\cdot h_n\cdot r \in S\cdot h_n$.
This shows $\pi_n \cdot R \subseteq S\cdot h_n$.

Step 2: $J \subseteq \pi_n \cdot R$: Take $a \in J$.
The fact that $a \in g_n\cdot S$ implies that $a=\pi_n \cdot s$ for some $s \in S$. By \zcref{lem:canc},
this implies that $s \in \mathcal{I}_m$ for all $m \ge 1$ such that $m \ne n$.
Here, \zcref{lem:canc} applies since $\pi_n\cdot s = a \in R$. Moreover, since $a = \pi_n\cdot s = v_n\cdot h_n\cdot s$ and $a \in S\cdot h_n$, we get $h_n\cdot s = v_n^{-1}\cdot a \in S\cdot h_n$, so $s \in \mathcal{I}_n$ by \zcref{lem:swap}. Hence $s \in \mathcal{I}_m$ for all $m \ge 1$, which gives $s \in R$, so $a = \pi_n\cdot s \in \pi_n\cdot R$.

Step 3: $R\cdot\pi_n = J$: By \zcref{prop:symmetry}, $\iota(R) = R$. Since $\iota$ reverses multiplication, and $\iota(g_n) = h_n$, $\iota(h_n) = g_n$, it maps $g_n\cdot S$ onto $S\cdot h_n$ and vice versa, so $\iota(J) = J$. Moreover,
\[
\iota(\pi_n) = \iota(v_n)\cdot\iota(g_n) = v_n\cdot h_n = \pi_n.
\]
Hence, by Steps 1 and 2, we obtain
$R\cdot\pi_n = \iota(\pi_n\cdot R) = \iota(J) = J$.
\end{proof}

The above will allow us to define a $\k$-algebra automorphism $R\to R$ which is a type of conjugation by a specific $\pi_n$:

\begin{lem}\label{lem:local-triv}
Let $ n \ge 1$.
\begin{enumerate}
\item For any $x \in R$, there is a unique $y \in R$ such that $x\cdot \pi_n = \pi_n \cdot y$. We will denote it by $y = c_n(x)$.
\item The map $c_n:R \to R$ defined in (1) is a $\k$-algebra automorphism.
\item Given any finite subset $F \subseteq \mathbb{Z}$ such that $n \notin F$,
if $x \in R \cap D_F[t]$, then $c_n(x) = x$.
\end{enumerate}
\end{lem}
\begin{proof}\quad
\begin{enumerate}
\item We have that $x \cdot \pi_n \in R\cdot \pi_n = \pi_n \cdot R$ by \zcref{lem:pi-left-right-ideals}.
Hence, there exists $y \in R$ such that $x\cdot \pi_n = \pi_n \cdot y$.
For uniqueness, if $\pi_n \cdot y_1 = \pi_n \cdot y_2$,
then $\pi_n\cdot (y_1 - y_2)=0$, and as $R$ is a domain and $\pi_n \ne 0$, we deduce that $y_1 = y_2$.
\item It is clear that $c_n$ is a $\k$-algebra homomorphism.
If $c_n(x) =0$, then $x\cdot \pi_n = 0$, and as $\pi_n \ne 0$ and $R$ is a domain, $x = 0$.
Finally, if $y \in R$, then $\pi_n \cdot y \in \pi_n \cdot R = R\cdot \pi_n$ by \zcref{lem:pi-left-right-ideals}.
Hence, there exists $x \in R$ such that $x\cdot \pi_n = \pi_n\cdot y$,
which implies that $c_n(x) =y$. Hence, $c_n$ is bijective.
\item By \zcref{lem:elempin}(2), we know that $x\cdot \pi_n = \pi_n \cdot x$, so $c_n(x) = x$ follows from the definition of $c_n$.
\end{enumerate}
\end{proof}

Next, we will compose different $c_n$'s to define the promised map $\tau:R \to R$. 
For any $N \ge 1$,
let 
\[
C^N := c_N\circ c_{N-1} \circ \dots \circ c_2 \circ c_1 :R \to R.
\]

\begin{lem}\label{lem:CNPN}
For all $N \ge 1$ and all $x \in R$,
it holds that $P_N \cdot C^N(x) = x\cdot P_N$,
where $P_N$ was defined in \zcref{eqn:PN}.
\end{lem}
\begin{proof}
We show this by induction on $N$.
If $N=1$, the claim is that $\pi_1 \cdot c_1(x) = x\cdot \pi_1$ which is the definition of $c_1$.
Assume the result holds for all $n < N$.
Then
\begin{align*}
P_N \cdot C^N(x) &= (\pi_1\cdot \dots \cdot \pi_{N-1}) \cdot \pi_N \cdot c_N(C^{N-1}(x)) \\ &= (\pi_1\cdot \dots \cdot \pi_{N-1}) \cdot C^{N-1}(x) \cdot \pi_N \\ 
&= x \cdot (\pi_1\cdot \dots \cdot \pi_{N-1}) \cdot \pi_N = x\cdot P_N.
\end{align*}
Here, the first equality is the definitions,
the second one is the defining property of $c_N$,
and the third is the induction hypothesis.
\end{proof}

We shall use this to prove an independence result:
\begin{lem}\label{lem:indepCN}
Let $x \in R$,
and let $F \subseteq \mathbb{Z}$ be a non-empty finite set such that $x \in D_F[t] \cap R$.
Let $N ,M \ge 1$ be such that $N \ge \max(F)$ and $M \ge \max(F)$.
Then there is an equality
$C^N(x) = C^M(x)$.
\end{lem}
\begin{proof}
We may assume that $M > N$, and take some $n > N$.
Notice that since by \zcref{lem:elempin}(3), the $\pi_n$'s commute with each other, it follows that $\pi_n$ commutes with $P_N$.
Since $n \notin F$, we also know by \zcref{lem:elempin}(2) that $\pi_n$ commutes with $x$.
This implies that:
\[
P_N \cdot C^N(x) \cdot \pi_n = x\cdot P_N \cdot \pi_n = x \cdot \pi_n \cdot P_N = \pi_n \cdot x \cdot P_N = \pi_n \cdot P_N \cdot C^N(x) = P_N \cdot \pi_n \cdot C^N(x),
\]
where we have used \zcref{lem:CNPN} twice. 
Cancelling the non-zero element $P_N$ on the left,
we deduce that 
\begin{equation}\label{eqn:CNpin}
C^N(x) \cdot \pi_n = \pi_n \cdot C^N(x)
\end{equation}
for all $n>N$.
This allows us to calculate:
\begin{align*}
P_M \cdot C^N(x) &= P_N \cdot (\pi_{N+1} \cdot \dots \cdot \pi_M) \cdot C^N(x) \\ &=
P_N \cdot C^N(x) \cdot  (\pi_{N+1} \cdot \dots \cdot \pi_M) \\ &= 
x\cdot P_N \cdot (\pi_{N+1} \cdot \dots \cdot \pi_M) = x\cdot P_M.
\end{align*}
Here, we used \zcref{eqn:CNpin} and \zcref{lem:CNPN}.
We thus see that
\[
P_M \cdot C^M(x) = x\cdot P_M = P_M \cdot C^N(x),
\]
and cancelling $P_M$ on the left this gives $C^M(x) = C^N(x)$.
\end{proof}

We now define $\tau:R \to R$ as follows.
For any $x \in R$, by \zcref{thm:divisionRing}(3) and the fact that $R$ is a subring of $S = D[t]$, there is a non-empty finite subset $F \subseteq \mathbb{Z}$, such that $x \in D_F[t] \cap R$.
Taking any $N > \max(F)$ with $N \ge 1$, we let:
\begin{equation}\label{eqn:tauDef}
\tau(x) = C^N(x) \in R.    
\end{equation}
This is independent of any choice.
Indeed, if $F_1, F_2 \subseteq \mathbb{Z}$ are two such non-empty sets such that $x \in D_{F_1}[t] \cap R$ and $x \in D_{F_2}[t] \cap R$,
and if $N_1 > \max(F_1)$, $N_1 \ge 1$ and $N_2 > \max(F_2)$, $N_2 \ge 1$, taking $M \ge \max(N_1,N_2)$, by \zcref{lem:indepCN}, we obtain:
\[
C^{N_1}(x) = C^M(x) = C^{N_2}(x),
\]
so we obtain a well defined function $\tau:R \to R$.

\begin{lem}\label{lem:tauinjective}
The map $\tau:R \to R$ defined in \zcref{eqn:tauDef} is an injective $\k$-algebra homomorphism.
\end{lem}
\begin{proof}
Suppose $x,y \in R$,
and take some non-empty finite subset $F \subseteq \mathbb{Z}$, such that $x,y \in D_F[t] \cap R$. Let $N \ge 1$ such that $N > \max(F)$.
Then $x\cdot y \in D_F[t]$,
so
\[
\tau(x\cdot y) = C^N(x \cdot y) = C^N(x) \cdot C^N(y) = \tau(x) \cdot \tau(y).
\]
Similarly, one shows that $\tau(x+y) = \tau(x) + \tau(y)$,
and that $\tau$ is $\k$-linear. Thus, $\tau$ is a $\k$-algebra homomorphism.
Finally, if $\tau(x) = 0$,
then $C^N(x) = 0$, but $C^N$, being a composition of injective maps,
is injective, so $x = 0$, and $\tau$ is injective.
\end{proof}

We need one more result about $c_n$ before we can prove that $\tau$ is an automorphism:

\begin{lem}\label{lem:cnDFinDF}
Let $F \subseteq \mathbb{Z}$ be non-empty finite set,
and let $n \in F$ be such that $n \ge 1$.
Then for any $x \in R\cap D_F[t]$,
it holds that $c_n(x) \in D_F[t]$ and $c_n^{-1}(x) \in D_F[t]$.
\end{lem}
\begin{proof}
Let us write $c_n(x) = \sum_{k=0}^m t^k \cdot a_k$.
Then:
\begin{equation}\label{eqn:pincnforcoeff}
\pi_n \cdot c_n(x) = (t\cdot v_n -u_n\cdot v_n) \cdot \sum_{k=0}^m t^k \cdot a_k = 
\sum_{k=0}^m (t^{k+1} \cdot v_n \cdot a_k - t^k \cdot u_n \cdot v_n \cdot a_k).
\end{equation}
On the other hand, by the definition of $c_n$,
we know that $\pi_n \cdot c_n(x) = x\cdot \pi_n$.
Since $n \in F$, we know that $\pi_n \in D_F[t]$.
Hence, $\pi_n \cdot c_n(x) = x\cdot \pi_n \in D_F[t]$.
Because of this, each of the coefficients of the polynomial in \zcref{eqn:pincnforcoeff} belongs to $D_F$.
In particular, the coefficient of $t^{m+1}$ in $x\cdot \pi_n$ is $v_n \cdot a_m \in D_F$,
so $a_m = v_n^{-1} \cdot (v_n \cdot a_m) \in D_F$.
Next, the coefficient of $t^m$ is $v_n\cdot a_{m-1} - u_n\cdot v_n \cdot a_m \in D_F$,
and we already know that $u_n\cdot v_n \cdot a_m \in D_F$, so we see that $v_n \cdot a_{m-1} \in D_F$, and hence, $a_{m-1} \in D_F$. Continuing in the same fashion,
we see that $a_0 ,\dots, a_m \in D_F$, which means that $c_n(x) \in D_F[t]$.
A nearly identical argument also shows the claim about $c_n^{-1}$.
\end{proof}

From this, we deduce:
\begin{lem}\label{lem:tauAuto}
The map $\tau:R \to R$ defined in \zcref{eqn:tauDef} is a $\k$-algebra automorphism.    
\end{lem}
\begin{proof}
In view of \zcref{lem:tauinjective}, it remains to show that $\tau$ is surjective.
Let $y \in R$, and let $F' \subseteq \mathbb{Z}$ be a finite set, such that $y \in D_{F'}[t]$.
Take some $N \ge 1$ such that $N \ge \max(F')$,
and modify $F'$, by letting $F = F' \cup \{1,\dots, N\}$.
Then, $N = \max(F)$, and $y \in D_F[t]$, by \zcref{thm:divisionRing}(3)(d).
Let $x = (c_1^{-1} \circ c_2^{-1} \circ \dots \circ c_N^{-1})(y)$, so $x \in R$.
Moreover, since $y \in R\cap D_F[t]$,
repeated use of \zcref{lem:cnDFinDF} tells us that $x \in D_F[t]$.
This implies by the definition of $\tau$ and the fact that $N = \max(F)$,
that
\[
\tau(x) = C^{N+1}(x) = c_{N+1}\circ c_N \circ \dots \circ c_1(x) = c_{N+1}(y).
\]
Since $N+1 \notin F$, and $y \in D_F[t]$,
by \zcref{lem:local-triv}(3), we get $c_{N+1}(y) =y$,
which shows that $\tau(x) = y$, and hence $\tau$ is surjective.
\end{proof} 

Here is the main result of this subsection:
\begin{prop}\label{prop:main-tausigma}
There exists a $\k$-algebra automorphism $\tau:R \to R$,
with the property that for all $r \in R$,
\[
\tau(\sigma(r)) \cdot \pi_0 = \pi_0 \cdot \sigma(\tau(r)),
\]
where $\pi_0$ was defined in \zcref{eqn:pi}.
\end{prop}
\begin{proof}
Let $\tau$ be the map defined in \zcref{eqn:tauDef}.
Then $\tau$ is a $\k$-linear automorphism by \zcref{lem:tauAuto}.
Let $r \in R$, and consider also $\sigma(r) \in R$.
Let $F \subseteq \mathbb{Z}$ be a finite set such that $r,\sigma(r) \in D_F[t] \cap R$,
and let $N \ge 1$ be such that $N > \max(F)$.
Then 
\[
\tau(r) = C^{N+1}(r), \quad \tau(\sigma(r)) = C^{N}(\sigma(r)).
\]
Hence, by \zcref{lem:CNPN},
we have equalities
\begin{equation}\label{eqn:ptaursigma}
P_{N+1} \cdot \tau(r) = r \cdot P_{N+1} ,\quad P_N \cdot \tau(\sigma(r)) = \sigma(r) \cdot P_N.
\end{equation}
According to \zcref{lem:elempin}(4), we know that $\beta(P_{N+1}) =\pi_0 \cdot P_N$.
As $\sigma$ is given by the restriction of $\beta$,
applying $\beta$ to $P_{N+1} \cdot \tau(r)$ we obtain:
\[
\sigma(P_{N+1} \cdot \tau(r)) =  \beta(P_{N+1} \cdot \tau(r)) = \beta(r\cdot P_{N+1}).
\]
So using the fact that $\beta$ respects multiplication and that it extends $\sigma$, this gives us:
\begin{equation}\label{eqn:beforesigtau}
\pi_0 \cdot P_N \cdot \sigma(\tau(r)) = \sigma(r) \cdot \pi_0 \cdot P_N.
\end{equation}
Going back to \zcref{eqn:ptaursigma}, we obtain from it that
\[
P_N \cdot \tau(\sigma(r)) \cdot \pi_0 = \sigma(r) \cdot P_N \cdot \pi_0.
\]
Since by \zcref{lem:elempin}(4), $\pi_0$ commutes with $P_N$, 
we obtain that
\[
P_N \cdot \tau(\sigma(r)) \cdot \pi_0 = \sigma(r) \cdot P_N \cdot \pi_0 = \sigma(r) \cdot \pi_0 \cdot P_N = \pi_0 \cdot P_N \cdot \sigma(\tau(r)) = P_N \cdot \pi_0 \cdot \sigma(\tau(r)),
\]
where the third equality is \zcref{eqn:beforesigtau}.
Cancelling the non-zero element $P_N$ from the left,
this shows that $\tau(\sigma(r)) \cdot \pi_0 = \pi_0 \cdot \sigma(\tau(r))$.
\end{proof}

\subsection{Dualizing complexes over trivial extensions}

Our next goal is to prove a general result about the existence of dualizing complexes over trivial extension rings over certain rings $R$. 

\subsubsection{Assumptions}\label{sec:assum} Throughout this subsection, we assume that $R$ is an associative ring,
which satisfies the following conditions:
\begin{enumerate}
\item $R$ is a $\k$-algebra.
\item $R$ is two-sided noetherian.
\item $R$ is two-sided Gorenstein, in the sense that $\injdim_R({}_R R) = \injdim_R(R_R) < \infty$.
\item There is an injective $\k$-algebra homomorphism $\sigma:R \to R$,
such that $1 \in (R:\sigma(R))_r\cdot R$,
and there is an element $q_0 \in R$ which is not a right zero divisor,
such that 
\[
\{r \in R \mid R\cdot r \subseteq \sigma(R)\} = R \cdot q_0.
\]
\item There is a $\k$-algebra automorphism $\tau:R \to R$,
such that $\tau(\sigma(r)) \cdot q_0 = q_0 \cdot \sigma(\tau(r))$ for all $r \in R$.
\end{enumerate}

Observe that by \zcref{prop:rightcondComp} and \zcref{prop:main-tausigma},
the ring $R$ from \zcref{thm:mainR} satisfies these conditions,
where $q_0 = \pi_0$. 
The element $q_0$ is not a right zero divisor because $R$ is also a domain,
but the domain assumption will no longer play a role in the sequel.
Also, $R$ is two-sided Gorenstein because it is hereditary.
Note that the point of this paper is that being Gorenstein is not left-right symmetric for noetherian rings, which is why we write ``two-sided Gorenstein'' explicitly.

Further note that being two-sided Gorenstein implies that the trivial bimodule ${}_R R_R$ is a dualizing complex over $R$. This follows easily from the definition of dualizing complexes.

First, we twist $R$ by $\tau$. 
We define $R_{\tau}$ to be the $R$-$R$-bimodule, which is equal to $R$ as a set,
and with left and right actions given by 
\[
x\circ r := x\cdot r, \quad   r\circ x := r\cdot\tau(x), \quad r \in R_{\tau}, x\in R.
\]
Thus, as a left $R$-module $R_{\tau}$ is simply ${}_R R$, but on the right the action twists by $\tau$. Note that since $\tau$ is $\k$-linear, $\k$ acts centrally on $R_{\tau}$.

In what follows, we will freely use the restriction functors $\restl, \restr$, introduced in \zcref{eqn:Restf}.

\begin{prop}\label{prop:rtaudualizing}
The bimodule $R_{\tau}$ is a dualizing complex over $R$.    
\end{prop}
\begin{proof}
The left $R$-module $\restl(R_{\tau})$ is just ${}_R R$,
which is clearly finitely generated and free,
and by assumption has finite injective dimension over $R$.
Thus, 
\[
\opn{Ext}^i_R(\restl(R_{\tau}),\restl(R_{\tau})) = \opn{Ext}^i_R(R,R) = 0
\]
for all $i>0$.
The homothety map $R^{\opn{op}} \to \opn{Hom}_R(\restl(R_{\tau}),\restl(R_{\tau}))$ is given by
$r \mapsto (a \mapsto a\circ r = a\cdot \tau(r))$.
Since $\tau$ is an isomorphism, this map is clearly bijective.
Next, we move to the right.
There, the homothety map $R \to \opn{Hom}_{R^{\opn{op}}}(\restr(R_{\tau}),\restr(R_{\tau}))$ is left multiplication, so clearly bijective.
The fact that $\tau$ is an isomorphism implies that the forgetful functor along $\tau$, denoted by
$\tau^*:\opn{Mod}_R \to \opn{Mod}_R$, is an equivalence of categories,
and 
\[
\tau^*(R_R) = \restr(R_{\tau}).
\]
Since $R_R$ is finitely generated projective and has finite injective dimension, and these properties are preserved by equivalences of module categories, we deduce that $\restr(R_{\tau})$ is also finitely generated projective and has finite injective dimension.
Thus, 
\[
\opn{Ext}^i_{R^{\opn{op}}}(\restr(R_{\tau}),\restr(R_{\tau})) = 0
\]
for all $i>0$, and $R_{\tau}$ is a dualizing complex over $R$.
\end{proof}

We let $A = R[x;\sigma]/(x^2) = R \skewtimes {}_{\sigma} R$ as in \zcref{eqn:hero}.
In particular, elements of $A$ are of the form $x+\epsilon\cdot y$ where $x,y \in R$, 
and $\epsilon^2 = 0$.
Our goal is to construct a dualizing complex over $A$.

Recall, as discussed in \zcref{sec:gen} that $\epsilon \cdot R \cong {}_{\sigma} R$ is naturally an $R$-$R$-bimodule. 
In the next lemma, the reader should think of $\omega$ as the dualizing complex over $R$,
and in that case, $K$ will eventually be the dualizing complex over $A$.

\begin{lem}\label{lem:bimoduleK}
Let $\omega$ and $N$ be two $R$-$R$-bimodules on which $\k$ acts centrally,
equipped with two $\k$-linear $R$-$R$-bimodule maps
\[
\lambda : \epsilon\cdot R \otimes_R N \to \omega, \quad \mu:N\otimes_R \epsilon \cdot R \to \omega.
\]
Then the following hold:
\begin{enumerate}
\item The $R$-$R$-bimodule $K = \omega \oplus N$, equipped with the $A$-actions:
\[
(x+\epsilon \cdot y) \circ (w,n) = (x\cdot w + \lambda(\epsilon \cdot y\otimes_R n), x\cdot n),
\]
and 
\[
(w,n) \circ (x+\epsilon \cdot y) = (w \cdot x + \mu(n \otimes_R \epsilon \cdot y), n\cdot x)
\]
is an $A$-$A$-bimodule, and $\k$ acts centrally on it.
\item The projection map $p:K\to \omega$ given by $p(w,n) = w$ is a $\k$-linear map of $R$-$R$-bimodules.
\end{enumerate}
\end{lem}
\begin{proof}
All these facts are routine calculations,
which follow from the $R$-linearity of $\lambda, \mu$ and the tensor product.
The fact that $\epsilon^2 =0$ means that no compatibility conditions between $\lambda, \mu$ are needed.
\end{proof}

Keeping in mind that $\omega$ will eventually be taken to be a dualizing complex over $R$,
our goal is to show that $K$ is both the left dual and the right dual of $A$ over $R$ with respect to $\omega$.

\begin{lem}\label{lem:twoduals}
Keep the notation of \zcref{lem:bimoduleK}.
\begin{enumerate}
\item
The map 
\[
\theta_{\ell} : K\to \opn{Hom}_R({}_R A_A, \omega) 
\]
given by
\[
\theta_{\ell}(k) (a) = p(a\circ k)
\]
is an $A$-$R$-bimodule homomorphism.
If moreover, the map 
\[
\lambda^* : N \to \opn{Hom}_R({}_R \epsilon \cdot R,\omega), \quad \lambda^*(n) (\epsilon \cdot y)  = \lambda(\epsilon\cdot y \otimes_R n)
\]
is bijective, then $\theta_{\ell}$ is also bijective.
\item
The map 
\[
\theta_{r} : K\to \opn{Hom}_{R^{\opn{op}}}({}_A A_R, \omega) 
\]
given by
\[
\theta_{r}(k) (a) = p(k\circ a)
\]
is an $R$-$A$-bimodule homomorphism.
If moreover, the map
\[
\mu^* :N \to \opn{Hom}_{R^{\opn{op}}}(\epsilon \cdot R_R,\omega), \quad \mu^*(n) (\epsilon \cdot y) = \mu(n\otimes_R \epsilon \cdot y)
\]
is bijective, then $\theta_r$ is also bijective.
\end{enumerate}
\end{lem}
\begin{proof}
We only prove (1), as the proof of the second claim is nearly identical.
First, we show that for $k \in K$, indeed $\theta_{\ell}(k)$ as defined above is an element of $\opn{Hom}_R({}_R A_A, \omega)$.
Take $r \in R$ and $a\in A$.
Since $p:K \to \omega$ is left $R$-linear, we have:
\[
\theta_{\ell}(k) (r\cdot a) = p((r\cdot a) \circ k) = p(r \circ (a \circ k)) = r\cdot p(a\circ k) = r\cdot \theta_{\ell}(k)(a).
\]
Next, we show that $\theta_{\ell}$ is a map of $A$-$R$-bimodules.
Take another element $a' \in A$.
Then 
\[
\theta_{\ell}(a' \circ k)(a) = p(a\circ a' \circ k) = \theta_{\ell}(k)(a\cdot a') = (a'\cdot \theta_{\ell}(k))(a),
\]
where the last equality is because of how the left $A$-structure on 
$\opn{Hom}_R({}_R A_A, \omega)$ is defined.
Now, take some $r \in R$.
Then
\[
\theta_{\ell}(k \circ r)(a) = p(a\circ k \circ r) = p(a\circ k)\cdot r,
\]
because $p$ is right $R$-linear.
Finally, suppose that $\lambda^*$ is bijective.
As left $R$-modules, we know that ${}_R A = R \oplus \epsilon \cdot R$.
This implies that the target of $\theta_{\ell}$ can be decomposed as
\[
\opn{Hom}_R({}_R A, \omega) = \omega \oplus \opn{Hom}_R(\epsilon\cdot R,\omega).
\]
Using this decomposition, 
an element $f \in \opn{Hom}_R({}_R A, \omega)$ decomposes as $(f(1),f|_{\epsilon \cdot R})$.
Then,
$\theta_{\ell}$ decomposes as
\[
\theta_{\ell} = 1_{\omega} \oplus \lambda^*.
\]
Indeed, for $k = (w,n) \in K$, we have
\[
\theta_{\ell}(w,n)(1) = p(1\circ (w,n)) = p(w,n) =w,
\]
while
\[
\theta_{\ell}(w,n)(\epsilon \cdot y) = p(\epsilon \cdot y \circ (w,n)) = 
p(\lambda(\epsilon \cdot y \otimes_R n),0) =\lambda(\epsilon \cdot y \otimes_R n) = \lambda^*(n)(\epsilon \cdot y).
\]
Using this presentation, it is clear that if $\lambda^*$ is bijective then so is $\theta_{\ell}$.
\end{proof}

A bimodule which is a dualizing complex is called a dualizing bimodule.
We now obtain a dualizing bimodule over $A$ given one over $R$:

\begin{lem}\label{lem:thedualcomp}
Keep the notation of \zcref{lem:bimoduleK} and \zcref{lem:twoduals},
and assume that the maps $\lambda^*$ and $\mu^*$ are bijective.
Suppose that $\omega$ is a dualizing bimodule over $R$.
Then $K$ is a dualizing bimodule over $A$,
and moreover,
\[
\injdim_A({}_A K) \le \injdim_R({}_R \omega), \quad
\injdim_A(K_A) \le \injdim_R(\omega_R).
\]
\end{lem}
\begin{proof}
First, we observe that $A$ is finitely generated projective over $R$ both on the left and on the right:
As left $R$-modules we simply have ${}_R A = R \oplus {}_{\sigma} R$,
and this was shown to be finitely generated projective in the proof of \zcref{prop:injdimleft}.
As right $R$-modules, we have that $A_R = R \oplus \epsilon \cdot R$ so $A$ is clearly finitely generated free of rank $2$.
The map $\theta_{\ell}$ gives us an isomorphism of left $A$-modules:
\[
{}_A K \cong \opn{Hom}_R({}_R A_A ,\omega).
\]
Assume that $\injdim_R({}_R \omega) = d$,
and let $J = 0 \to I^0 \to \dots \to I^d \to 0$ be an injective resolution of $\omega$ as a left module over $R$.
Applying the exact functor $\opn{Hom}_R({}_R A_A,-)$ to the quasi-isomorphism $\omega \to J$, 
we obtain a sequence of quasi-isomorphisms of left $A$-modules:
\[
{}_A K \cong \opn{Hom}_R({}_R A_A ,\omega) \cong \opn{Hom}_R({}_R A_A,J).
\]
Since $\opn{Hom}_R({}_R A_A,J)$ is a complex of injective left $A$-modules concentrated in cohomological degrees $[0,d]$,
we deduce that $\injdim_A({}_A K) \le d$.
The inequality $\injdim_A(K_A) \le \injdim_R(\omega_R)$ on the right is proven similarly using the fact that $\theta_r$ is an isomorphism.
Next, we show that $K$ is finitely generated over $A$.
The isomorphism $\theta_r : K \to \opn{Hom}_{R^{\opn{op}}}({}_A A_R,\omega)$ tells us that, as a left $R$-module, $K$ is a direct summand of some $\omega^n$, because $A$ itself, being finitely generated projective, is a direct summand of some $R^n$.
Thus, $K$ is finitely generated as a left module over $R$,
so a fortiori it is finitely generated over the larger ring $A$.
The proof that it is finitely generated on the right is symmetric.
It remains to compute $\mrm{R}\opn{Hom}_A(K,K)$ (as well as over $A^{\opn{op}}$).
First, we compute the higher $\opn{Ext}$'s.
For $i > 0$,
using the isomorphisms $\theta$ and adjunction (which is valid since $A$ is projective over $R$),
we obtain
\[
\opn{Ext}^i_A(K,K) \cong \opn{Ext}^i_A(K,\opn{Hom}_R(A,\omega)) \cong
\opn{Ext}^i_R(K,\omega),
\]
but we have seen that $K$ is a direct summand of $\omega^n$ as a left $R$-module,
and $\opn{Ext}^i_R(\omega,\omega) = 0$ since $\omega$ is a dualizing bimodule.
Thus, $\opn{Ext}^i_A(K,K) = 0$ for all $i>0$, and similarly on the right.
It remains to show that the homothety maps are bijective.
We first compute using the maps $\theta$:
\[
\opn{Hom}_A(K,K) \cong \opn{Hom}_A(K,\opn{Hom}_R(A,\omega)) \cong 
\opn{Hom}_R(K,\omega) \cong \opn{Hom}_R(\opn{Hom}_{R^{\opn{op}}}(A_R,\omega),\omega).
\]
Using this identification,
the homothety map $A^{\opn{op}} \to \opn{Hom}_A(K,K)$ is isomorphic to the map
\[
A^{\opn{op}} \to \opn{Hom}_R(\opn{Hom}_{R^{\opn{op}}}(A_R,\omega),\omega).
\]
Indeed, under the identification $\opn{Hom}_A(K,K) \cong \opn{Hom}_R(K,\omega)$, $\varphi \mapsto p\circ\varphi$, the homothety of $a \in A$ is sent to the map $k \mapsto p(k\circ a) = \theta_r(k)(a)$. Identifying $K$ with $\opn{Hom}_{R^{\opn{op}}}(A_R,\omega)$ via $\theta_r$, it becomes the evaluation map $f \mapsto f(a)$.

But $A$ is projective, 
so $\opn{Hom}_{R^{\opn{op}}}(A_R,\omega) = \mrm{R}\opn{Hom}_{R^{\opn{op}}}(A_R,\omega)$,
while
\[
\opn{Hom}_R(\opn{Hom}_{R^{\opn{op}}}(A_R,\omega),\omega) = \mrm{R}\opn{Hom}_R(\opn{Hom}_{R^{\opn{op}}}(A_R,\omega),\omega),
\]
because we have seen that the higher $\opn{Ext}$'s vanish.
Thus, the obtained map is just the duality map from the right $R$-module $A_R$ to its derived double dual with respect to the dualizing bimodule $\omega$, and hence, it is bijective (\cite[Proposition 1.3]{YZ1}). The proof for the other homothety map is identical.
\end{proof}

All that remains is to construct an $R$-bimodule $N$ and maps $\lambda, \mu$ which satisfy the above conditions.
We let $N := R_{\tau\circ \sigma}$.
Thus, $N$ is the $R$-$R$-bimodule which is equal to $R$ as an abelian group, and with actions given by:
\[
r \circ n = r\cdot n, \quad n \circ r = n\cdot \tau(\sigma(r)), \quad r\in R, n\in N.
\]
Let $\omega = R_{\tau}$. By \zcref{prop:rtaudualizing},
$\omega$ is a dualizing bimodule over $R$.

\begin{lem}\label{lem:mapmu}
Define a map $\mu:N\otimes_R \epsilon\cdot R \to R_{\tau}$ by $\mu(n\otimes_R \epsilon \cdot y) = n \cdot \tau(y)$.
Then $\mu$ is a well defined $\k$-linear $R$-$R$-bimodule map,
and the induced map
\[
\mu^* :N \to \opn{Hom}_{R^{\opn{op}}}(\epsilon \cdot R_R,\omega), \quad \mu^*(n) (\epsilon \cdot y) = \mu(n\otimes_R \epsilon \cdot y)
\]
is bijective.
\end{lem}
\begin{proof}
The map $\mu$ is well defined as the corresponding map $N\times \epsilon \cdot R \to R_{\tau}$ is $R$-biadditive.
Indeed, note that for $n\otimes_R \epsilon \cdot y \in N\otimes_R \epsilon\cdot R$ and $r \in R$,
we have:
\[
\mu(n\circ r \otimes_R \epsilon\cdot y) = n\cdot\tau(\sigma(r))\cdot\tau(y) = n\cdot\tau(\sigma(r)\cdot y) = \mu(n\otimes_R \epsilon\cdot\sigma(r)\cdot y) = \mu(n\otimes_R r\cdot\epsilon\cdot y),
\]
where the last equality is by \zcref{eqn:multA}.

Since $\tau$ is $\k$-linear, so is $\mu$.
Let us check that $\mu$ is an $R$-$R$-bimodule map.
First, on the left:
Take $r \in R$ and $n\otimes_R \epsilon \cdot y \in N\otimes_R \epsilon\cdot R$.
Then:
\[
\mu(r\circ (n\otimes_R \epsilon \cdot y)) = 
r\cdot n \cdot \tau(y) = r \circ \mu(n \otimes_R \epsilon \cdot y).
\]
More interesting is the right side, where the twist occurs.
There, taking $r \in R$ and $n\otimes_R \epsilon \cdot y \in N\otimes_R \epsilon\cdot R$ we obtain:
\[
\mu((n\otimes_R \epsilon \cdot y) \circ r) = 
n\cdot \tau(y\cdot r) = n\cdot \tau(y)\cdot \tau(r) = (n\cdot \tau(y)) \circ r = \mu(n\otimes_R \epsilon \cdot y) \circ r.
\]
It remains to show that $\mu^*$ is bijective.
Note that $\epsilon \cdot R$ is a free right $R$-module of rank $1$,
which implies that the map 
$\opn{Hom}_{R^{\opn{op}}}(\epsilon \cdot R_R,\omega) \to \omega$ given by
$f \mapsto f(\epsilon)$ is bijective.
Composing $\mu^*$ with this bijection, we get
$\mu^*(n)(\epsilon) = \mu(n\otimes_R \epsilon) = n\cdot \tau(1) = n$
which is just the identity on the set $R$.
Hence, $\mu^*$ is also bijective.
\end{proof}

Next, we define $\lambda$. 
The element $q_0 \in R$ was introduced in \zcref{sec:assum}(4).

\begin{lem}\label{lem:maplambda}
Define a map $\lambda:\epsilon\cdot R\otimes_R N \to R_{\tau}$ by
$\lambda(\epsilon \cdot y\otimes_R n) = \sigma^{-1}(y\cdot n \cdot q_0)$.
Then $\lambda$ is a well defined $\k$-linear $R$-$R$-bimodule map, and the induced map 
\[
\lambda^* : N \to \opn{Hom}_R({}_R \epsilon \cdot R,\omega), \quad \lambda^*(n) (\epsilon \cdot y)  = \lambda(\epsilon\cdot y \otimes_R n)
\]
is bijective.
\end{lem}
\begin{proof}
First, we must explain why writing $\sigma^{-1}(y\cdot n \cdot q_0)$ makes sense.
This is because 
\[
y\cdot n \cdot q_0 \in R \cdot q_0 = \{r \in R \mid R\cdot r \subseteq \sigma(R)\} \subseteq \sigma(R).
\]
Hence, there exists $z \in R$ such that $\sigma(z) = y\cdot n \cdot q_0$,
and since $\sigma$ is injective, that $z$ is unique, and $z = \sigma^{-1}(y\cdot n \cdot q_0)$.
Clearly $\lambda$ is additive in both tensor arguments. 
We check that $\lambda$ is $R$-biadditive so it makes sense to define such a map on a tensor product:
for $r \in R$, $\epsilon \cdot y \in \epsilon \cdot R$ and $n \in N$,
we have:
\[
\lambda(\epsilon\cdot y \cdot r \otimes_R n) = \sigma^{-1}(y\cdot r \cdot n \cdot q_0) = \lambda(\epsilon\cdot y \otimes_R r \circ n).
\]
Next, we check that $\lambda$ is a map of $R$-$R$-bimodules.
First, on the left:
\begin{align*}
\lambda(r\cdot (\epsilon\cdot y)\otimes_R n) &= 
\lambda(\epsilon\cdot \sigma(r)\cdot y\otimes_R n) \\ &= 
\sigma^{-1} (\sigma(r)\cdot y\cdot n \cdot q_0) = r\cdot \sigma^{-1}(y\cdot n \cdot q_0) = r\cdot \lambda(\epsilon\cdot y\otimes_R n).
\end{align*}
Now, on the right:
\begin{align*}
& \lambda(\epsilon\cdot y\otimes_R n\circ r) = \sigma^{-1}(y\cdot n \circ r \cdot q_0) \\ &= \sigma^{-1}(y\cdot n\cdot \tau(\sigma(r)) \cdot q_0) = 
\sigma^{-1}(y\cdot n \cdot q_0 \cdot \sigma(\tau(r))) \\ &=
\sigma^{-1}(y\cdot n \cdot q_0) \cdot \tau(r) =
\lambda(\epsilon\cdot y\otimes_R n )\circ r.
\end{align*}
Here, the second equality is the definition of the right action on $N$,
the third one is \zcref{sec:assum}(5),
and the fifth is the definition of the right action on $R_{\tau}$.
Since $\sigma, \tau$ are $\k$-linear, we see that $\lambda$ is a well defined $\k$-linear $R$-$R$-bimodule map,
and all that remains is to show that $\lambda^*$ is bijective.
First, $\lambda^*$ is injective:
if $\lambda^*(n) = 0$, in particular 
\[
0 = \lambda^*(n)(\epsilon \cdot 1) = \lambda(\epsilon\cdot 1\otimes_R n) = \sigma^{-1}(n\cdot q_0).
\]
Applying $\sigma$, this implies that $n\cdot q_0 = 0$,
but $q_0$ is not a right zero-divisor, so $n=0$ and $\lambda^*$ is injective.
To show $\lambda^*$ is surjective, take some $f \in \opn{Hom}_R({}_R \epsilon \cdot R,\omega)$.
Define a map $g:R \to R$ by $g(y) = \sigma(f(\epsilon \cdot y))$.
Considering $R$ as a left $\sigma(R)$-module, note that $g$ is a map of left $\sigma(R)$-modules, because for $r,y \in R$, we have:
\begin{align*}
g(\sigma(r)\cdot y) &= \sigma(f(\epsilon\cdot \sigma(r)\cdot y)) = 
\sigma(f(r\cdot \epsilon \cdot y)) \\
& = \sigma(r\cdot f(\epsilon \cdot y)) = \sigma(r) \cdot \sigma(f(\epsilon \cdot y) ) = \sigma(r) \cdot g(y).
\end{align*}
Set $q := g(1) = \sigma(f(\epsilon))$.
Since $g$ is $\sigma(R)$-linear, 
we may apply \zcref{lem:fgProj}(2), and deduce that 
\[
g(y) = y\cdot g(1) = y\cdot q
\]
for all $y \in R$.
Since $\opn{Im}(g) \subseteq \sigma(R)$, we obtain that $R\cdot q \subseteq \sigma(R)$, so by \zcref{sec:assum}(4),
there exists $n\in R$ such that $q = n\cdot q_0$.
Then, we may calculate:
\[
f(\epsilon \cdot y) = \sigma^{-1}(g(y)) = \sigma^{-1}(y\cdot q) =
\sigma^{-1}(y\cdot n \cdot q_0) = \lambda(\epsilon\cdot y\otimes_R n),
\]
which shows that $\lambda^*(n) = f$, so $\lambda^*$ is bijective.
\end{proof}

Here is the main result of this section:

\begin{thm}\label{thm:dualcomplex}
Let $\k$ be any field with $\opn{char}(\k) \ne 2$,
and let $A$ be the $\k$-algebra constructed in \zcref{cor:main}.
Then $A$ has a dualizing complex $K$.
Moreover, $K$ satisfies:
\begin{enumerate}
\item $K$ is a single $A$-$A$-bimodule concentrated in degree $0$, a dualizing bimodule.
\item $\injdim_A({}_A K) \le 1$, and $\injdim_A(K_A) \le 1$.
\item The field $\k$ acts on $K$ centrally.
\end{enumerate}
\end{thm}
\begin{proof}
We observed in \zcref{sec:assum} that the ring $R$ from \zcref{thm:mainR} satisfies all the assumptions of \zcref{sec:assum}.
Thus, the results of this subsection apply to $R$.
Since $R$ is hereditary, $\injdim_R({}_R R) \le 1$, $\injdim_R(R_R) \le 1$.
Let $\omega = R_{\tau}$.
By \zcref{prop:rtaudualizing}, it is a dualizing complex over $R$,
and since the forgetful functor $\tau^*$ is an equivalence of categories that sends $R_R$ to $\omega_R$, 
we deduce that $\injdim_R(\omega_R) \le 1$.
As a left $R$-module, ${}_R \omega = {}_R R$, so we already know that $\injdim_R({}_R \omega) \le 1$.
Applying \zcref{lem:thedualcomp}, which holds because of \zcref{lem:mapmu} and \zcref{lem:maplambda}, 
we obtain that the $A$-$A$-bimodule $K$ constructed in \zcref{lem:bimoduleK} is a dualizing bimodule over $A$ which satisfies all the conclusions of this theorem.
\end{proof}

\section{Failure of Rickard-Keller conditions}\label{sec:RK}

Let $\mathcal{T}$ be a triangulated category which is closed under coproducts and products. In our applications, $\mathcal{T}$ will always be $\cat{D}({}_B \opn{Mod})$ or $\cat{D}(\opn{Mod}_B)$ for some ring $B$.

Given a non-empty set of objects $\mathbf{S} \subseteq \mathcal{T}$,
the localizing subcategory of $\mathcal{T}$ generated by $\mathbf{S}$,
denoted by $\opn{Loc}_{\mathcal{T}}(\mathbf{S})$ is the smallest triangulated subcategory of $\mathcal{T}$ which contains $\mathbf{S}$ and is closed under arbitrary coproducts.
If $\opn{Loc}_{\mathcal{T}}(\mathbf{S}) = \mathcal{T}$, we say that $\mathbf{S}$ generates $\mathcal{T}$.
Dually, replacing coproducts by products, 
the smallest triangulated subcategory of $\mathcal{T}$ which contains $\mathbf{S}$ and is closed under arbitrary products is denoted by $\opn{Coloc}_{\mathcal{T}}(\mathbf{S})$, and when it is equal to $\mathcal{T}$, we say that $\mathbf{S}$ cogenerates $\mathcal{T}$.

For a ring $B$, we say that injectives generate for $B$ (on the left) if the set ${}_B \opn{Inj}$ of injective left $B$-modules generates $\cat{D}({}_B \opn{Mod})$. Similarly, we say that projectives cogenerate for $B$ if the set ${}_B \opn{Proj}$ of projective left $B$-modules cogenerates
$\cat{D}({}_B \opn{Mod})$. 
Injective cogeneration and projective generation are defined similarly,
and hold for any ring (see \cite[Proposition 2.2]{Rickard} for the projective generation result. Injective cogeneration is proved similarly).

Any commutative noetherian ring satisfies injective generation \cite[Theorem 3.3]{Rickard}, 
a fact that follows from the Hopkins-Neeman \cite{Hopkinsglobal,Neeman92} classification of localizing subcategories of $\cat{D}({}_B \opn{Mod})$ when $B$ is commutative noetherian.
Moreover, by \cite[Theorem 3.2]{Rickard}, if $B$ is any ring such that $\injdim_B({}_B B) < \infty$, then injectives generate for $B$.

In 2001, Keller \cite{Keller01} observed that for a finite-dimensional algebra, injective generation implies several of the classical homological conjectures, and asked if this is the case for the finitistic dimension conjecture. This was answered positively by Rickard in \cite[Theorem 4.3]{Rickard}, and was later generalized to artin algebras having a dualizing complex in \cite{Zhang}, 
and, more recently, \cite{ShFinDim}, by the author, to all two-sided noetherian rings with a dualizing complex.
Thus, we now know that for a two-sided noetherian ring $B$, which has a dualizing complex,
if injectives generate for $B$ (on the left), then the big finitistic projective dimension of $B$
\[
\opn{FPD}(B) := \sup_{M \in {}_B\opn{Mod}, \projdim_B(M) < \infty}(\projdim_B(M))
\]
is finite. 
Similarly, Rickard showed in \cite[Proposition 5.2]{Rickard} that for finite-dimensional algebras, 
projective cogeneration on the right implies finite finitistic dimension on the left.
And that result was also generalized by the author to cover all two-sided noetherian rings with a dualizing complex, \cite[Theorem 5.1]{ShGor}. Rings with finite injective dimension and finite finitistic dimension on the opposite side are known to have finite injective dimension on both sides by \cite[Proposition 6.10]{AusRei}.

By the Rickard-Keller conditions we mean injective generation and projective cogeneration.
In \cite[Theorem 6.3]{ShGor}, using the above, we gave necessary and sufficient conditions for a two-sided noetherian ring with a dualizing complex having finite injective dimension on the left to have finite injective dimension on the right. Since this is violated for the algebra $A$ constructed in this paper, we obtain:

\begin{thm}\label{thm:failure}
Let $\k$ be any field with $\opn{char}(\k) \ne 2$,
and let $A$ be the $\k$-algebra constructed in \zcref{cor:main}.
Then the following additional properties hold for $A$:
\begin{enumerate}
\item Injectives generate for $A$ on the left.
\item Injectives do not generate for $A$ on the right.
That is, $\opn{Loc}_{\cat{D}(\opn{Mod}_A)}(\opn{Inj}_A) \ne \cat{D}(\opn{Mod}_A)$.
\item $\opn{Coloc}_{\cat{D}({}_A \opn{Mod})}({}_A \opn{Proj}) \ne \cat{D}({}_A \opn{Mod})$. 
That is, projectives do not cogenerate for $A$ on the left.
\item $\opn{Coloc}_{\cat{D}({}_A \opn{Mod})}({}_A \opn{Flat}) \ne \cat{D}({}_A \opn{Mod})$. 
That is, flats do not cogenerate for $A$ on the left.
\item There exists a bounded above cochain complex $J$ of injective left $A$-modules which is acyclic but not contractible.
\item There exists a bounded below cochain complex $P$ of projective right $A$-modules which is acyclic but not contractible.
\item There exists an injective left $A$-module $I$ with infinite flat dimension.
\item There exists a dualizing complex $\omega'$ over $A$ which has infinite flat dimension as a complex of left $A$-modules.
\end{enumerate}
\end{thm}
\begin{proof}
Since by \zcref{cor:main}, we know that 
\[
\injdim_A({}_A A) = 1, \quad \injdim_A(A_A) = +\infty,
\]
and since by \zcref{thm:dualcomplex}, $A$ has a dualizing complex,
all conditions except for (1), (2) and (6) follow immediately from \cite[Theorem 6.3]{ShGor}.
(1) follows from \cite[Theorem 3.2]{Rickard},
and (5) is equivalent to (6) by \cite[Theorem 5.1]{ShGor}.
Finally, (2) follows from (3) also by \cite[Theorem 5.1]{ShGor}.
\end{proof}

\begin{rem}
The definition of a dualizing complex used in \cite{ShFinDim, ShGor} is slightly stronger than the one used in this paper. This is because the proofs there use the covariant Grothendieck duality \cite{InKr,Neeflat}, which requires moreover that a dualizing complex $\omega$ over a ring $B$ is isomorphic to a complex of $B$-$B$-bimodules, such that simultaneously 
$\restlb(\omega)$ is a complex of injective left $B$-modules,
and $\restrb(\omega)$ is a complex of injective right $B$-modules.
However, when one assumes in addition that $B$ is a $\k$-algebra over a field (as is the case for our algebra $A$), replacing $\omega$ by an injective resolution over the enveloping algebra $B\otimes_{\k} B^{\opn{op}}$ gives the required stronger variant of a dualizing complex,
so we may invoke the above results for $A$ freely.
\end{rem}

\begin{rem}
Regarding conditions (5) and (6) of \zcref{thm:failure}, 
one should mention that by \cite[Theorem 5.1]{ShGor},
the condition that any bounded above acyclic complex of injective left $B$-modules is contractible is equivalent to $\opn{FPD}(B^{\opn{op}}) < \infty$ for any two-sided noetherian ring $B$ with a dualizing complex. 
Equivalently, any bounded below acyclic complex of projective right $B$-modules is contractible.
These are fairly new homological conditions, first introduced in \cite{ShGor}.
In the appendix to \cite{ShGor},
Nakamura and Thompson showed \cite[Corollary A.2 and Corollary A.8]{ShGor} these conditions hold over any commutative noetherian ring.
Examples of complexes violating these conditions were constructed in \cite{Pos2},
but all of them were over noncommutative non-noetherian rings.
This result shows that not only do such violations exist over two-sided noetherian rings,
but they can even exist when the injective dimension is finite on the left.
\end{rem}

\section{The Auslander condition}\label{sec:Auslander}

A two-sided noetherian ring $B$ is said to satisfy the Auslander condition if for every finitely generated left $B$-module $M$,
for all $i \ge 0$,
and 
any right $B$-module $N$ which is a submodule of the right $B$-module $\opn{Ext}^i_B(M,B)$, it holds that $\opn{Ext}^j_{B^{\opn{op}}}(N,B) = 0$ for all $0 \le j < i$. Other equivalent formulations of this condition are given in \cite[Theorem 3.7]{TExtBook},
where it is also shown that this condition holds for $B$ if and only if it holds for $B^{\opn{op}}$.
One of the equivalent conditions,
which is the one we shall use,
says that in the minimal injective resolution $0\to B \to J^0 \to J^1 \to \dots$ of $B$ as a left module over itself,
it holds that $\flatdim_B({}_B J^i) \le i$ for all $i\ge 0$.

Auslander and Reiten showed in \cite[Corollary 5.5]{AuReiGor} that if an artin algebra $B$ has 
\[
\injdim_B({}_B B) <\infty,
\]
and if $B$ satisfies the Auslander condition, then $B$ satisfies Gorenstein symmetry,
i.e., 
\[
\injdim_B({}_B B) = \injdim_B(B_B) < \infty.
\]
This was extended by Huang in \cite[Proposition 4.6]{Huang} to all two-sided artinian rings,
and even more recently by Shen in \cite[Theorem 1.3]{Shen26}, to the much larger class of all two-sided noetherian rings which are module-finite algebras over a central commutative noetherian ring.
Huang explicitly asked \cite[p.~2]{Huang} if this implication holds for any two-sided noetherian ring.

The purpose of this section is to show that the algebra $A$ of \zcref{cor:main} satisfies the Auslander condition.
In particular, for general two-sided noetherian rings,
Auslander's condition \textbf{does not imply} Gorenstein symmetry,
and Huang's question has a negative answer.

To do that, we will only need to make the assumptions of \zcref{thm:main-gen-const}. Thus, we make the following assumptions:
Let $R$ be a two-sided noetherian domain which is left hereditary, and is not a division ring. 
Let $\sigma:R \to R$ be an injective ring homomorphism which is not surjective,
and assume that 
$1 \in (R:\sigma(R))_r\cdot R$.
Let $A = R[x;\sigma]/(x^2)$ as in \zcref{eqn:hero}.
We will show that, under these assumptions, $A$ satisfies the Auslander condition.

Since $R$ is a two-sided noetherian domain, it is a left and right Ore domain \cite[Theorem 2.1.15]{NCBook}, so it has a division ring of fractions, in which every element is both a right and a left fraction of elements of $R$, see \cite[Corollaries 2.1.14 and 2.1.4]{NCBook}.
We denote its division ring of fractions by $K$.

The condition that $1 \in (R:\sigma(R))_r\cdot R$
implies there exists some $0 \ne q_0 \in (R:\sigma(R))_r$.
We will fix, from now on, such an element $0 \ne q_0 \in R$.
Thus, it satisfies that 
\begin{equation}\label{eqn:qocondAS}
r\cdot q_0 \in \sigma(R) \quad \mbox{for all $r \in R$.}
\end{equation}
In particular, taking $r=1$, we see that $q_0 \in \sigma(R)$.
Let $p_0 \in R$ be the unique element such that $\sigma(p_0) = q_0$.
Since $\sigma$ is injective and $q_0 \ne 0$, also $p_0 \ne 0$.

\begin{lem}
The map $\sigma:R \to R$ extends uniquely to a ring map $\sigma:K \to K$.
Moreover, this extension is an automorphism of $K$.
\end{lem}
\begin{proof}
Composing $\sigma:R \to R$ with the inclusion $R \inj K$,
we get a map $\sigma:R \to K$.
As $\sigma:R \to K$ is injective, it maps nonzero elements of $R$ to unit elements of $K$.
By \cite[Lemma 2.1.4]{NCBook}, the universal property of Ore localization,
it extends uniquely to a ring homomorphism $K \to K$, which we also denote by $\sigma$.
Since $K$ is a division ring, $\sigma:K \to K$ is still injective. 
It remains to show it is surjective.
For any $r \in R$, in the division ring $K$ there is an equality
$r = (r\cdot q_0) \cdot q_0^{-1}$.
By \zcref{eqn:qocondAS}, $r\cdot q_0 \in \sigma(R)$,
so there is $s \in R$ with $\sigma(s) = r\cdot q_0$.
Thus, in $K$ we have 
\[
\sigma(s\cdot p_0^{-1}) = (r\cdot q_0) \cdot q_0^{-1} = r.
\]
Finally, if $x \in K$ is arbitrary,
we may write $x = x_1 \cdot x_2^{-1}$ with $x_1,x_2 \in R$, and $x_2 \ne 0$,
and we have just seen there exist $y_1,y_2 \in K$ such that $\sigma(y_1) = x_1$, $\sigma(y_2) = x_2$, and $x_2 \ne 0$ implies $y_2 \ne 0$.
Thus, $\sigma(y_1\cdot y_2^{-1}) = x_1\cdot x_2^{-1} = x$,
so $\sigma$ is surjective.
\end{proof}

Let ${}_{\sigma} K$ be the $K$-$K$-bimodule whose underlying abelian group is equal to $K$,
and with left and right actions given by
\[
r\circ k := \sigma(r) \cdot k, \quad k \circ r := k \cdot r, \quad r\in K, k \in {}_{\sigma} K.
\]
Using it, we construct the trivial extension ring $Q := K \skewtimes {}_{\sigma} K$. As in \zcref{sec:gen}, its elements are of the form $x + \epsilon \cdot y$, with $x,y \in K$, and multiplication defined with the same formula as in \zcref{eqn:multA}. In particular, 
since the multiplication formulas are the same, $R \inj K$,
and $\sigma:K \to K$ extends $\sigma:R \to R$,
and as $A$ was defined to be $R \skewtimes {}_{\sigma} R$,
it follows that
there is a natural inclusion $A \inj Q$.
We will show that it is the first term in the injective resolution of $A$.
For this, first we show flatness:

\begin{lem}\label{lem:QisFlat}
The inclusion $A \inj Q$ is a left and right Ore localization 
with respect to the set $T := \{x+\epsilon \cdot y\mid x,y \in R, x \ne 0\}$.
In particular, $Q$ is flat over $A$, both as a left $A$-module and as a right $A$-module.
\end{lem}
\begin{proof}
Denote by $\varphi:A \inj Q$ the inclusion map.
According to \cite[Definition 2.1.3]{NCBook}, 
to show $\varphi$ is a right Ore localization with respect to $T$,
we must check that $\varphi$ maps any element of $T$ to a unit element of $Q$, that any element of $Q$ can be written as $\varphi(r)\cdot \varphi(s)^{-1}$ for some $r\in A$, $s \in T$,
and that $\ker(\varphi) = \{a\in A \mid \exists s \in T, a\cdot s =0\}$.
We will show this, as well as the left version of these conditions.
If $x+\epsilon \cdot y \in T$, so $x \ne 0$,
note that
\[
(x+\epsilon \cdot y) \cdot (x^{-1}-\epsilon \cdot \sigma(x)^{-1}\cdot y \cdot x^{-1}) = 1 + \epsilon(-\sigma(x)\cdot \sigma(x)^{-1}\cdot y \cdot x^{-1} + y \cdot x^{-1}) = 1.
\]
Similarly, multiplication in reverse also gives $1$,
which shows that $\varphi$ maps $T$ to units in $Q$.
Take $x+\epsilon \cdot y \in Q$, so $x,y \in K$.
According to \cite[Proposition 2.1.16(i)]{NCBook},
we may write $x=a\cdot s^{-1}$, $y = b\cdot s^{-1}$ for some $a,b, s\in R$, with $s \ne 0$.
Thus, $(a+\epsilon \cdot b)\cdot s^{-1} = x + \epsilon\cdot y$,
and clearly $s \in T$, proving the claim.
Next, we show the same claim on the left for $x+\epsilon \cdot y$, $x,y \in K$.
Let us write $x = s_1^{-1} \cdot a_1$ and $\sigma^{-1}(y) = s_2^{-1} \cdot b_2$ for some $a_1, b_2, s_1 ,s_2 \in R$ where $s_1 \ne 0$ and $s_2 \ne 0$. Using the fact that $R\setminus \{0\}$ is a left Ore set in $R$, we may find non-zero elements $r_1, r_2 \in R$, such that $r_1\cdot s_1 = r_2 \cdot s_2$. Denote this nonzero element of $R$ by $s$.
Finally, let $a = r_1 \cdot a_1$ and $b = \sigma(r_2\cdot b_2)$.
Then $s^{-1}\cdot a = s_1^{-1} \cdot r_1^{-1} \cdot r_1 \cdot a_1 = x$,
while $\sigma(s)^{-1} \cdot b = \sigma(r_2\cdot s_2)^{-1}\cdot \sigma(r_2\cdot b_2) = \sigma(s_2)^{-1}\cdot \sigma(b_2) = y$.
Hence, using the definition of left multiplication, as in \zcref{eqn:multA}, we see:
\[
s^{-1}\cdot (a+\epsilon \cdot b) = s^{-1}\cdot a + \epsilon \cdot \sigma(s)^{-1} \cdot b = x+\epsilon \cdot y.
\]
Finally, $\varphi$ being an inclusion is injective, so its kernel is $\{0\}$. On the other hand, for $a = x_1+\epsilon\cdot y_1\in A$, if $a\cdot s = 0$ for some $s \in T$, say $s = x_2+\epsilon \cdot y_2$, with $x_2 \ne 0$,
then
\[
0 = (x_1 + \epsilon y_1) \cdot (x_2 + \epsilon y_2) = x_1\cdot x_2 + \epsilon(\sigma(x_1)\cdot y_2+ y_1\cdot x_2).
\]
Hence, $x_1\cdot x_2 = 0$, but $R$ is a domain and $x_2 \ne 0$, so $x_1 = 0$. Moreover, 
\[
0 = \sigma(x_1)\cdot y_2+ y_1\cdot x_2 = y_1\cdot x_2,
\]
so the same argument shows that $y_1 =0$, and $a = 0$.
On the left, if $s \cdot a =0$ where $a = x_1+\epsilon\cdot y_1\in A$,
and $s = x_2 + \epsilon \cdot y_2 \in T$ (so $x_2 \ne 0$),
then
\[
0 = (x_2 + \epsilon y_2) \cdot (x_1 + \epsilon y_1) = x_2\cdot x_1 + \epsilon(\sigma(x_2)\cdot y_1+ y_2\cdot x_1).
\]
Hence, $x_2 \cdot x_1 = 0$, so $x_1 = 0$.
This implies $\sigma(x_2)\cdot y_1 = 0$, and as $\sigma$ is injective, $\sigma(x_2) \ne 0$, so $y_1 = 0$, proving that $\varphi$ is a left and right Ore localization with respect to $T$.
The flatness claim is \cite[Proposition 2.1.16(ii)]{NCBook} (applied once on each side).
\end{proof}

Next, we show self-injectivity, which together with flatness will immediately imply injectivity:

\begin{lem}\label{lem:QisInj}
The regular left module ${}_Q Q$ is injective over $Q$.
In other words, $Q$ is a self-injective ring on the left.
\end{lem}
\begin{proof}
We shall use Baer's criterion to test injectivity. To do that, let us first find all left ideals in $Q$.
Let $J$ be such left ideal,
and let $x + \epsilon \cdot y \in J$ be some non-zero element.
Then $x,y \in K$, and if $x \ne 0$,
the computation in the proof of \zcref{lem:QisFlat} shows that such an element is a unit in $Q$, with inverse being $(x^{-1}-\epsilon \cdot \sigma(x)^{-1}\cdot y \cdot x^{-1})$ (there, it was further assumed that $x,y \in R$, but this is not needed for this computation to hold).
Thus, any such $J$ is equal to $Q$.
Otherwise, we see that $J \subseteq \epsilon\cdot K$.
Assuming $J \ne 0$, there is some $0 \ne y \in K$ such that $\epsilon \cdot y \in J$, and then 
\[
Q\cdot \epsilon \cdot y = \{\epsilon \cdot \sigma(a) \cdot y \mid a \in K\} = \epsilon \cdot K,
\]
where the last equality follows because $\sigma$ is surjective and $K$ is a division ring.
Thus, we see that $\epsilon \cdot K$ is the only non-trivial left ideal in $Q$. To use Baer's criterion, let $f:\epsilon \cdot K \to Q$ be a left $Q$-linear map, and we must extend it to a map $Q \to Q$.
Let $e = f(\epsilon)$, and write $e = e_1 + \epsilon \cdot e_2$, $e_1 ,e_2 \in K$. Then $0 = f(\epsilon\cdot \epsilon) = \epsilon \cdot e = \epsilon \cdot e_1$. Hence, $e_1 = 0$. Thus, $e = \epsilon \cdot e_2$.
For any $\epsilon \cdot x \in \epsilon\cdot K$, with $x \in K$, let $y = \sigma^{-1}(x) \in K$. This makes sense since $\sigma$ is bijective.
Then $y\cdot \epsilon = \epsilon \cdot \sigma(y) = \epsilon \cdot x$.
We have that $f(\epsilon \cdot x) = f(y\cdot \epsilon) = y \cdot e = y \cdot \epsilon \cdot e_2 = \epsilon \cdot x \cdot e_2$.
This shows that $f$ may be extended to the map $Q \to Q$ given by $z\mapsto z\cdot e_2$, which is clearly left $Q$-linear. Hence, ${}_Q Q$ is injective over $Q$.
\end{proof}

Since by \zcref{lem:QisFlat}, the right $A$-module $Q_{A}$ is flat, and since by \zcref{lem:QisInj}, 
${}_Q Q$ is injective,
the adjunction
\[
\opn{Hom}_A(-,{}_A Q) \cong \opn{Hom}_Q(Q_{A}\otimes_A -,{}_Q Q)
\]
shows that ${}_A Q$ is injective over $A$.

Recall that an inclusion of left $A$-modules $N \inj M$ is an essential extension if the intersection of $N$ with any non-zero submodule of $M$ is non-zero.

\begin{lem}
The inclusion of left $A$-modules ${}_A A \inj {}_A Q$ is an essential extension.    
\end{lem}
\begin{proof}
Let $L \subseteq {}_A Q$ be a non-zero left submodule,
and let $0 \ne x \in L$.
By \zcref{lem:QisFlat}, we may write $x = s^{-1}\cdot a$ with $a \in A$, $s \in T$. Then $a = s\cdot x \in L$, and $x \ne 0$ implies $a \ne 0$,
so $0 \ne a \in L \cap A$.
\end{proof}

Here is the main result of this section:
\begin{thm}\label{thm:Auslander}
Let $R$ be a two-sided noetherian domain which is left hereditary, and is not a division ring. 
Let $\sigma:R \to R$ be an injective ring homomorphism which is not surjective,
and assume that 
$1 \in (R:\sigma(R))_r\cdot R$.
Let $A = R[x;\sigma]/(x^2)$ as in \zcref{eqn:hero}.
Then $A$ satisfies the Auslander condition.
\end{thm}
\begin{proof}
Since ${}_A Q$ is an essential extension of $A$ and is injective over $A$, we deduce that ${}_A Q$ is an injective hull of ${}_A A$.
Consider the short exact sequence of left $A$-modules:
\begin{equation}\label{eqn:injResA}
0 \to {}_A A \inj {}_A Q \to {}_A Q/A \to 0.
\end{equation}
By \zcref{thm:main-gen-const}, we know that $\injdim_A({}_A A) =1$,
so the fact that $Q$ is an injective left $A$-module containing ${}_A A$ implies that 
the first cosyzygy of \zcref{eqn:injResA}, $\opn{coker}(A \inj Q) = {}_A Q/A$ is injective over $A$.
Since ${}_A Q$ is an injective hull of ${}_A A$,
we deduce that \zcref{eqn:injResA} is the minimal injective resolution of ${}_A A$.
To verify the Auslander condition, it remains to show that $\flatdim_A ({}_A Q) = 0$, which is \zcref{lem:QisFlat}, and that $\flatdim_A ({}_A Q/A) \le 1$. The latter follows directly from exactness of \zcref{eqn:injResA} and the fact that ${}_A A$ and ${}_A Q$ are flat over $A$.
\end{proof}

We immediately obtain:
\begin{cor}\label{cor:Auslander}
Let $\k$ be any field with $\opn{char}(\k) \ne 2$,
and let $A$ be the $\k$-algebra constructed in \zcref{cor:main}.
Then $A$ satisfies the Auslander condition.    
In particular, Huang's question has a negative answer.
\end{cor}
\begin{proof}
By \zcref{thm:mainR}, the ring $R$ satisfies the assumptions of \zcref{thm:Auslander}.
\end{proof}

\textbf{Statement on the use of AI.} The construction used in \zcref{thm:mainR} was obtained from interactions of the author with Claude Fable 5.1 and Claude Opus 5.5.
The paper itself was fully human written.

\bibliographystyle{abbrv}
\bibliography{main}

\end{document}